\documentclass[11pt,reqno]{article}
\usepackage{latexsym, amsmath, amssymb, a4, epsfig}
\usepackage{graphicx}
\usepackage{color}

\newtheorem{theorem}{Theorem}[section]

\newtheorem{lemma}[theorem]{Lemma}
\newtheorem{prop}[theorem]{Proposition}

\newcommand{\qed}{ $\Box$}
\newcommand{\ds}{\displaystyle}
\newcommand{\pf}{\noindent {\sl Proof}. \ }
\newcommand{\p}{\partial}
\newcommand{\pd}[2]{\frac {\p #1}{\p #2}}
\newcommand{\eqnref}[1]{(\ref {#1})}

\newcommand{\Cbb}{\mathbb{C}}
\newcommand{\Ebb}{\mathbb{E}}

\newcommand{\Kbb}{\mathbb{K}}
\newcommand{\Rbb}{\mathbb{R}}
\newcommand{\Sbb}{\mathbb{S}}

\newcommand{\la}{\langle}
\newcommand{\ra}{\rangle}

\newcommand{\Ccal}{\mathcal{C}}
\newcommand{\Dcal}{\mathcal{D}}
\newcommand{\Ecal}{\mathcal{E}}
\newcommand{\Fcal}{\mathcal{F}}
\newcommand{\Hcal}{\mathcal{H}}
\newcommand{\Ical}{\mathcal{I}}

\newcommand{\Kcal}{\mathcal{K}}
\newcommand{\Mcal}{\mathcal{M}}
\newcommand{\Pcal}{\mathcal{P}}

\newcommand{\Scal}{\mathcal{S}}

\def\Bc{{\bf c}}
\def\Bd{{\bf d}}

\def\Bf{{\bf f}}
\def\Bg{{\bf g}}

\newcommand{\Ga}{\alpha}

\newcommand{\Gd}{\delta}
\newcommand{\Ge}{\epsilon}

\newcommand{\Gvf}{\varphi}

\newcommand{\Gl}{\lambda}

\newcommand{\Gt}{\theta}

\newcommand{\Gs}{\sigma}

\newcommand{\Gx}{\xi}

\newcommand{\Gz}{\zeta}
\newcommand{\GD}{\Delta}

\newcommand{\GG}{\Gamma}
\newcommand{\GL}{\Lambda}

\newcommand{\GO}{\Omega}

\newcommand{\beq}{\begin{equation}}
\newcommand{\eeq}{\end{equation}}

\numberwithin{equation}{section}
\numberwithin{figure}{section}

\begin{document}
\title{Spectral resolution of the Neumann-Poincar\'{e} operator on intersecting disks and analysis of plasmon resonance\thanks{\footnotesize This work is supported by the Korean Ministry of Education, Sciences and Technology through NRF grants Nos. 2010-0017532 (to H.K) and 2012003224 (to S.Y) and by the Korean Ministry of Science, ICT and Future Planning through NRF grant No. NRF-2013R1A1A3012931 (to M.L).}}

\author{Hyeonbae Kang\thanks{Department of Mathematics, Inha University, Incheon
402-751, S. Korea (hbkang@inha.ac.kr).} \and Mikyoung Lim\thanks{Department of Mathematical Sciences,
Korea Advanced Institute of Science and Technology, Daejeon 305-701, Korea (mklim@kaist.ac.kr, shyu@kaist.ac.kr).}\and Sanghyeon Yu\footnotemark[3]}

\maketitle

\begin{abstract}
The purpose of this paper is to investigate the spectral nature of the Neumann-Poincar\'e operator on the intersecting disks, which is a domain with the Lipschitz boundary. The complete spectral resolution of the operator is derived, which shows in particular that it admits only the absolutely continuous spectrum, no singularly continuous spectrum and no pure point spectrum. We then quantitatively analyze using the spectral resolution the plasmon resonance at the absolutely continuous spectrum.
\end{abstract}



%%%%%%%%%%%%%%%%%%%%%%%%%%%%%%%%%%%%%%%%%%%%%
\section{Introduction}\label{sec:intro}
%%%%%%%%%%%%%%%%%%%%%%%%%%%%%%%%%%%%%%%%%%%%%

The Neumann-Poincar\'e (NP) operator is a boundary integral operator naturally arising when solving classical Dirichlet or Neumann boundary value problems using layer potentials. This operator is not self-adjoint with respect to the usual $L^2$ inner product unless the domain is a disk or a ball. However, it is recently found that by introducing a new inner product the NP operator can be symmetrized and admits a spectral resolution. If the boundary of the domain is smooth, then the corresponding NP operator is compact, and hence has only the point spectrum accumulating to $0$. But, if the boundary of the domain is Lipschitz, then the NP operator is not compact (it is a singular integral operator), and it is quite interesting to investigate the nature of the spectrum of the NP operator, in particular, whether it admits both point and continuous spectrum. In this paper we derive a complete spectral resolution of the NP operator on the intersecting disks which is a Lipschitz domain, from which we are able to show that the NP operator on the intersecting disks has only absolutely continuous spectrum. As far as we are aware of, this is the first example of a domain with the Lipschitz boundary for which the complete spectrum of the NP operator is derived.

The study of the NP operator goes back to Poincar\'e as the name of the operator alludes. It has been a central object in the theory of singular integral operators developed during the last century (see \cite{kenig}). Recently there is growing interest in the spectral property of the NP operator in relation to plasmonics. The plasmonic structure consists of an inclusion of the negative dielectric constant (with dissipation) embedded in the matrix of the positive dielectric constant. The eigenvalue of the NP operator is related to the plasmon resonance (see for example \cite{Grieser}). The negativity of the dielectric constant of the inclusion makes it necessary to look into the spectrum of the corresponding NP operator. If the dielectric constant is positive, which corresponds to the ellipticity of the problem, only the resolvent of the NP operator matters. In this paper we derive quantitatively precise estimates of the plasmon resonance at the continuous spectrum of the NP operator on intersecting disks in the quasi-static limit. Particularly we obtain precise rates of resonance as the dissipation parameter tends to zero. As far as we are aware of, this is the first such a result. We mention that the plasmon resonance (with a non-zero frequency) at the continuous spectrum on the intersecting disks was also studied in \cite{pendry}. But there quantitative estimate of resonance is missing.

As mentioned before the NP operator is not self-adjoint with respect to the usual $L^2$ inner product, unless the domain is a disc or a ball (see \cite{Li01}). However, it is recently found in \cite{KhPuSh07} (see also \cite{Kang}) that NP operator can be symmetrized using Plemelj's symmetrization principle. There are a few domains with smooth boundaries on which the spectrum of the NP operator is known; a complete spectrum is known on disks, ellipses, and balls, and a few eigenvalues on ellipsoids (see \cite{Kang} and references therein). The spectrum of the NP operator on two disjoint disks has been computed recently, and applied to analysis of high concentration of the stress in between closely located disks \cite{BT, BT2, LY1}. The spectral property of the NP operator is barely known. But, the bounds on the essential spectrum of the NP operator on the curvilinear polygonal domains have been derived in \cite{PP}. The intersecting disks is among domains considered in the paper, and interestingly the result of this paper reveals that the bounds there are optimal for the intersecting disks. Recently, the spectral theory of the NP operator finds new applications. In \cite{ACKLM} it was efficiently used for analysis of the cloaking by anomalous localized resonance on coated plasmonic structure. More recently, uniformity (with respect to the elliptic constants) of elliptic estimates and boundary perturbation formula has been proved using the spectral property of the NP operator \cite{KKLS}. It is also applied to analysis of plasmon resonance on smooth domains \cite{Ando}.

The finding of this paper shows that the NP operator on intersecting disks does not admits a point spectrum on $H_0^{-1/2}(\p\GO)$. We suspect that this is rather exceptional and domains with corners, especially 4 or more corners like squares, admit point spectrums. This is an interesting question to investigate. In this regards, it is also interesting to see whether all the non-smooth Lipschitz domains admit non-empty continuous spectrum. It is worth mentioning that the NP operator on $H^{-1/2}(\p\GO)$ space on any domain $\GO$ (with possibly multiple connected components) with the Lipschitz boundary has $1/2$ as an eigenvalue and its multiplicity is the same as the number of connected boundary components \cite{ACKLY}.

This paper is organized as follows. In the next section we introduce the NP operator and the bipolar coordinate system which is an essential ingredient of this paper. In section \ref{sec3} we explicitly compute the single layer potential and the NP operator on the intersecting disks in terms of the bipolar coordinates. In section \ref{sec4} we derive a resolution of the identity and the spectral resolution of the NP operator. The last two sections are to derive estimates for the plasmon resonance.

%%%%%%%%%%%%%%%%%%%%%%%%%%%%%%%%%%%%%%%%%%%%%
\section{Preliminaries}
%%%%%%%%%%%%%%%%%%%%%%%%%%%%%%%%%%%%%%%%%%%%%

Let us fix some notation first. Let $\GO$ be a bounded domain in $\Rbb^2$ with the Lipschitz boundary. We denote by $H^{-1/2}(\p{\GO})$ the dual space of $H^{1/2}(\p{\GO})$ where the latter denotes the usual Sobolev space. We use the notation $\la \cdot , \cdot \ra$ for either the $L^2$ inner product or the duality  pairing of $H^{-1/2}$ and $H^{1/2}$ on $\p\GO$. We denote by $\| \cdot \|_s$ the Sobolev norm on $H^s(\p\GO)$ for $s=-1/2$ or $1/2$. Let $H_0^{-1/2}(\p{\GO})$ be the space of $\psi \in H^{-1/2}(\p{\GO})$ satisfying $\la \psi, 1 \ra =0$.

%%%%%%%%%%%%%%%%%%%%%%%%%%%%%%%%%%%%%%%%%%%%%%%%
\subsection{The NP operator and symmetrization}
%%%%%%%%%%%%%%%%%%%%%%%%%%%%%%%%%%%%%%%%%%%%%%%%

For $x=(x_1,x_2) \neq 0$, let
$\GG(x)$ be the fundamental solution to the
Laplacian, {\it i.e.},
\beq
\GG (x) =\frac{1}{2\pi} \ln |x|.
\eeq
A classical way of solving the Neumann boundary value problem on $\GO$ is to use the single layer potential $\Scal_{\p\GO} [\Gvf]$ which is defined by
\beq\label{slp}
\Scal_{\p\GO} [\Gvf] (x) := \int_{\p \GO} \Gamma (x-y) \Gvf (y) \, d\sigma(y) \;, \quad x \in \Rbb^2 ,
\eeq
where $\Gvf \in L^2(\p\GO)$ or $H^{-1/2}(\p\GO)$. It is well known (see, for example, \cite{AmKa07Book2, Fo95}) that $\Scal_{\p\GO} [\Gvf]$ satisfies the jump relation
\beq\label{singlejump}
\frac{\p}{\p\nu} \Scal_{\p\GO} [\Gvf] \big |_\pm (x) = \biggl( \pm \frac{1}{2} I + \Kcal_{\p\GO}^* \biggr) [\Gvf] (x),
\quad \mbox{a.e. } x \in \p\GO\;,
\eeq
where the operator $\Kcal_{\p\GO}^*$ is defined by
\beq \label{introkd2}
\Kcal^*_{\p\GO} [\Gvf] (x) = \frac{1}{2\pi} \int_{\p \GO}
\frac{(x -y) \cdot \nu_x}{|x-y|^2} \Gvf(y)\,d\sigma(y)\;.
\eeq
Here $\pm$ indicates the limits (to $\p\GO$) from outside and inside of $\GO$, respectively, and $\nu_x$ denotes the outward unit normal vector to $\p\GO$ at $x \in \p\GO$. The operator $\Kcal_{\p\GO}^*$ is called the Neumann-Poincar\'e (NP) operator associated with the domain $\GO$. Sometimes it is called the adjoint NP operator to distinguish it from the NP operator $\Kcal_{\p\GO}$ ($\Kcal_{\p\GO}$ is the $L^2$ adjoint of $\Kcal_{\p\GO}^*$).

Even though $\Kcal_{\p\GO}^*$ is not self-adjoint on the usual $L^2$-space, it can be symmetrized by introducing a new inner product. Recently it is proved in \cite{KhPuSh07} that $\Kcal_{\p\GO}^*$ can be symmetrized using Plemelj's symmetrization principle (also known as Calder\'on's identity)
\beq\label{Plemelj}
\Scal_{\p \GO} \Kcal^*_{\p \GO} = \Kcal_{\p \GO} \Scal_{\p \GO}.
\eeq
In fact, if we define
\beq\label{innerp}
\la \Gvf,\psi \ra_{\Hcal^*}:= -\la \Gvf,  \Scal_{\p \GO}[\psi] \ra
\eeq
for $\Gvf, \psi \in H^{-1/2}_0(\p \GO)$, then it is an inner product on $H^{-1/2}_0(\p \GO)$ and $\Kcal^*_{\p \GO}$ is self-adjoint with respect to this inner product.

Let $\Hcal^*=\Hcal^*(\p \GO)$ be the space $H^{-1/2}_0(\p \GO)$ equipped with the inner product $\la \ , \ \ra_{\Hcal^*}$ and denote the norm associated with
$\la \ , \ \ra_{\Hcal^*}$ by $\| \, \cdot\,\|_{\Hcal^*}$. It is proved in \cite{KKLS} that $\| \, \cdot\,\|_{\Hcal^*}$ is equivalent to $H^{-1/2}$ norm, {\it i.e.},
there are constants $C_1$ and $C_2$ such that
\beq\label{Hcal-1/23}
C_1 \| \Gvf \|_{-1/2} \le  \| \Gvf \|_{\Hcal^*}  \le C_2\| \Gvf \|_{-1/2}.
\eeq
for all $\Gvf  \in H^{-1/2}_0(\p \GO)$.

Since $\Kcal_{\p \GO}^*$ is self-adjoint on $\Hcal^*$, its spectrum $\Gs(\Kcal^*_{\p \GO})$ is real and consists of point and continuous spectra. Moreover, by the spectral resolution theorem there is a family of projection operators $\Ecal_t$ on $\Hcal^*$ (called a resolution of the identity) such that
\beq\label{specresol1}
\Kcal^*_{\p \GO} = \int_{-1/2}^{1/2} t \, d\Ecal_t.
\eeq
See \cite{Yo}.  We emphasize that the spectrum of $\Kcal_{\p\GO}^*$ on $\Hcal^*$ lies in $(-\frac{1}{2}, \frac{1}{2})$ (see \cite{Ke29}).

%%%%%%%%%%%%%%%%%%%%%%%%%%%%%%%%%%%%%%%%%%%%%%
\subsection{Intersecting disks and bipolar coordinates}
%%%%%%%%%%%%%%%%%%%%%%%%%%%%%%%%%%%%%%%%%%%%%%

\begin{figure}[!ht]
\begin{center}
\epsfig{figure=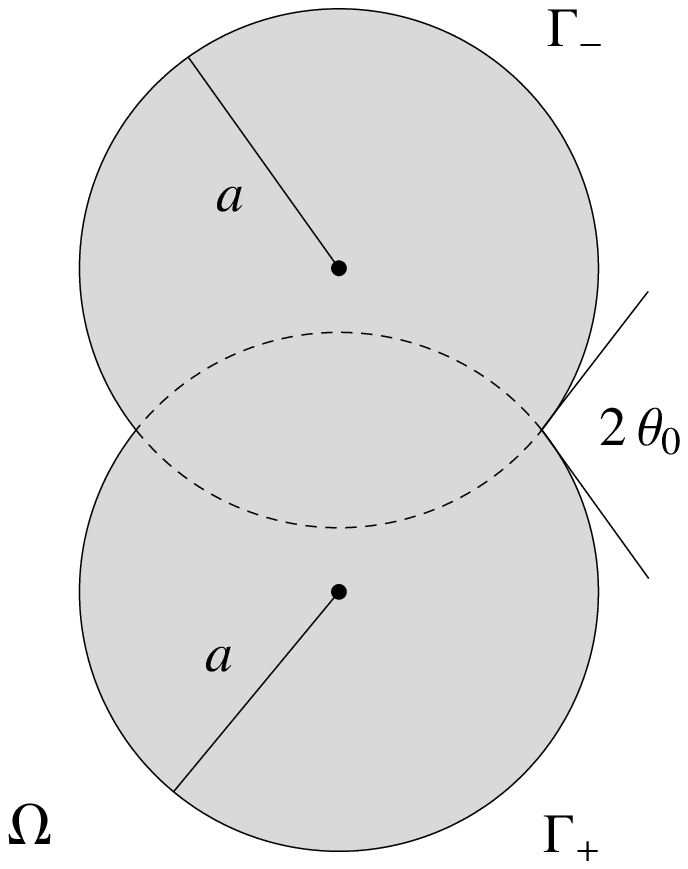,width=5.8cm}\hskip 1cm\epsfig{figure=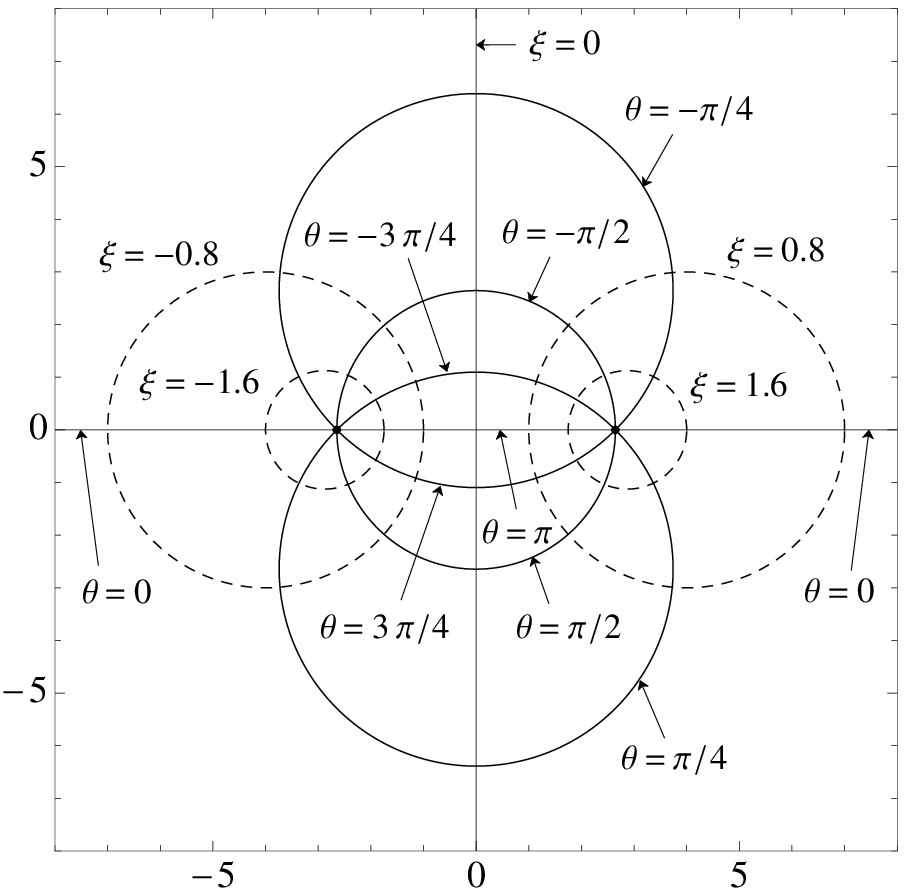,width=6cm}
\end{center}
\caption{The left figure is the intersecting disks $\GO$. Level coordinate curves of $\Gz$ (when $\alpha=\sqrt7$) are illustrated in $z$-plane in the right figure.}\label{intersecting}
\end{figure}

Let $B_1=B_{a}(\Bc_1)$ and $B_2=B_{a}(\Bc_2)$ be two disks of the same radii $a$ centered at $\Bc_1$ and $\Bc_2$, respectively, which are intersecting, namely,  $|\Bc_1-\Bc_2| < 2a$. The domain $\GO$ we consider in this paper is defined to be
\beq
\GO:=B_1\cup B_2.
\eeq
(We may consider $B_1 \cap B_2$ instead.) Note that $\p\GO$ is Lipschitz continuous with corners at the intersecting points. We assume that $\Bc_1$ is on the $x_2$-axis and $\Bc_1=-\Bc_2$ so that the intersecting points lie on the $x_1$-axis and they are symmetric with respect to the origin.
Let $2 \Gt_0$ ($\Gt_0 >0$) be the angle between two circles at the intersecting point (see Fig.\;\ref{intersecting}).
One can easily see that the intersecting points are given by $\pm a \sin \Gt_0$.

\begin{figure}[!ht]
\begin{center}
\epsfig{figure=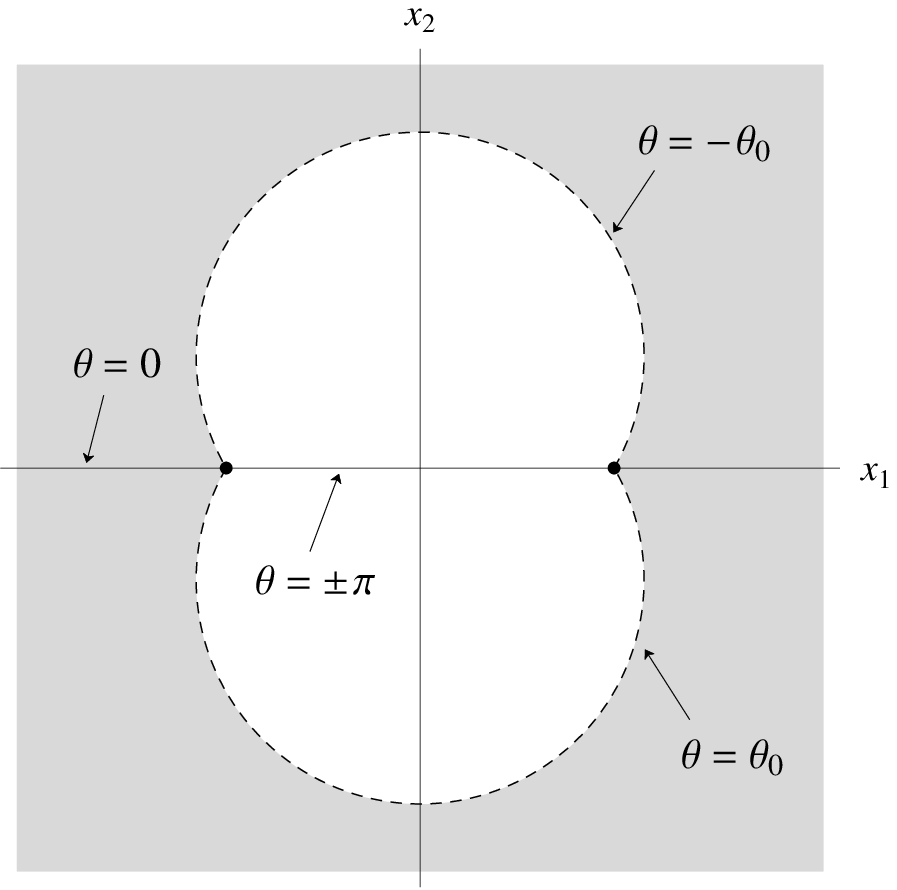,width=5.8cm} \hskip  1cm
\epsfig{figure=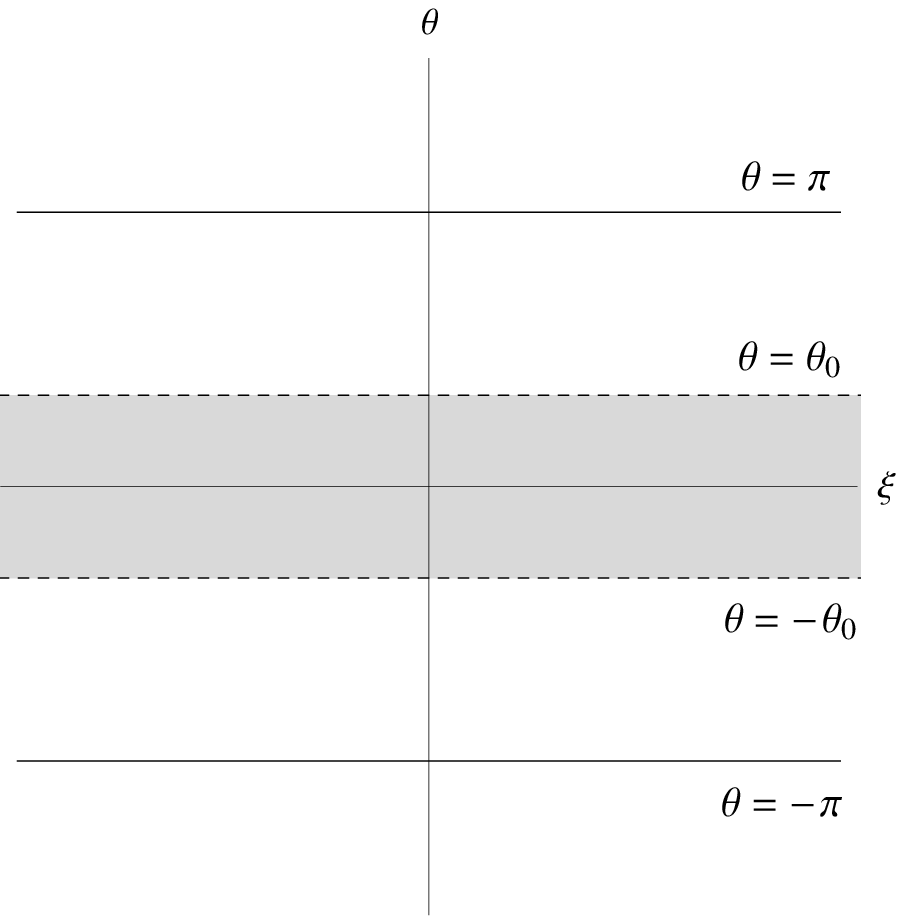,width=5.8cm}
\end{center}
\caption{The exterior of the intersecting disks (the gray region in the left figure) is the image of the infinite strip (the gray region in the right figure) via the conformal map $z=x_1+ix_2=\Phi(\xi+i\theta)$.} \label{conformal}
\end{figure}

For a given positive number $\Ga$ the bilinear transformation
\beq
\Gl(z) = \frac{z+\Ga}{z-\Ga}
\eeq
is a conformal map from the Riemann sphere onto itself, and maps the interval $[-\Ga, \Ga]$ onto the negative real axis. So,
\beq\label{loglz}
\Gz = \xi+ i \Gt:= \mbox{Log}\, \Gl(z)
\eeq
is well defined for $z \in \Cbb \setminus [-\Ga, \Ga]$, where $\mbox{Log}$ is the logarithm with the principal branch. The coordinate $\Gz$ or $(\xi, \Gt)$ is the bipolar coordinates we will use. The level coordinate curves for $\Gz$ are illustrated in Fig.\;\ref{intersecting}.
Here and afterwards, we take
\beq
\Ga= a \sin \Gt_0,
\eeq
so that $(\Ga,0)$ and $(-\Ga,0)$ are two intersecting points of the circles. Note that $\mbox{Log} \Gl$ maps $\Cbb \setminus [-\Ga, \Ga]$ onto $\{ \Gl=\Gx+i\Gt, \, -\infty< \xi < \infty, -\pi < \Gt < \pi \}$. If we let $\GG_+$ and $\GG_-$ be the lower and upper halves of $\p\GO$, respectively, one can see that $\mbox{Log} \Gl$ maps $\GG_+$ onto $\{\Gt=\Gt_0\}$, and $\GG_-$ onto $\{\Gt=-\Gt_0\}$, and outside of $\GO$ onto $\{ \Gl=\Gx+i\Gt, \, -\infty< \xi < \infty, -\Gt_0 < \Gt < \Gt_0 \}$ (see Fig.\;\ref{conformal}). Note that we denote the lower half by $\GG_+$ because it is mapped to $\{\Gt=+\Gt_0\}$.

Let $\Phi$ be the inverse of $\mbox{Log} \Gl$, {\it i.e.},
\beq\label{def_Phi}
z=x_1+i x_2= \Phi(\Gz):=\alpha \frac{e^\Gz+1}{e^\Gz-1}.
\eeq
We also regard $\Phi$ as a mapping from $\Rbb^2$ into $\Rbb^2$ by identifying $x=(x_1,x_2)\in \Rbb^2$ with $z=x_1+ix_2 \in \Cbb$.
Here we recall a few facts to be used later. We define the scale factor $h$ by
\beq\label{def_h}
h(\xi,\Gt):=\frac{1}{|\Phi'(\xi+i\Gt)|}=\frac{\cosh\xi-\cos\Gt}{\alpha}.
\eeq
Since $\Phi$ is a conformal transformation, $h$ is the Jacobian of the coordinate change. In particular, one can see that the outward normal derivative on $\p\GO$ is transformed in the following way:
\beq\label{nor_bipolar}
\pd{u}{\nu} = \mp h(\xi,\Gt_0)\pd{(u \circ \Phi)}{\theta} \quad\mbox{on } \GG_{\pm},
\eeq
and the line element $d\Gs$ on $\p\GO$ is transformed by
\beq\label{jacobian}
d\Gs = \frac{1}{h(\xi,\Gt_0)} d\xi.
\eeq

Define
\beq\label{Psidef}
\mbox{Log} \Gl(z) =:\Psi_1(z) + i\Psi_2(z).
\eeq
Then we have
\beq\label{inverse}
\begin{cases}
\Phi(\Psi_1(x), \Gt_0)=x \quad &\mbox{if } x \in \GG_+, \\
\Phi(\Psi_1(x), -\Gt_0)=x \quad &\mbox{if } x \in \GG_-.
\end{cases}
\eeq

%%%%%%%%%%%%%%%%%%%%%%%%%%%%%%%%%%%%%%%%%%%%%%%%%%%%%%%%%%%%%%%%%%
\section{Some explicit computations}\label{sec3}
%%%%%%%%%%%%%%%%%%%%%%%%%%%%%%%%%%%%%%%%%%%%%%%%%%%%%%%%%%%%%%%%%%

%%%%%%%%%%%%%%%%%%%%%%%%%%%%%%%%%%%%%%%%%%%%%%%%%%%%%%%%%%%%%%%%%%
\subsection{Single layer potential on $\Rbb^2$}
%%%%%%%%%%%%%%%%%%%%%%%%%%%%%%%%%%%%%%%%%%%%%%%%%%%%%%%%%%%%%%%%%%

The purpose of this section is to derive an explicit formula for the single layer potential on $\p\GO$ using the bipolar coordinates.
For $\Gvf \in \Hcal^*$ let
\beq\notag
u(x):= \Scal_{\p\GO} [\Gvf](x), \quad x \in \Rbb^2.
\eeq
Then $u$ is the solution to the following transmission problem:
\beq\label{transmission_eigen}
\begin{cases}
\ds\Delta u = 0 \quad& \mbox{in } \GO \cup (\mathbb{R}^2\setminus \overline{\GO}),\\
\ds u|_+=u|_- \quad& \mbox{on }  \p\GO, \\
\ds \frac{\p u}{\p \nu}\Big|_+-\frac{\p u}{\p \nu}\Big|_-=\varphi \quad&\mbox{on } \p\GO,\\
\ds u(x)=O(|x|^{-1}) \quad& \mbox{as } |x| \rightarrow \infty.
\end{cases}
\eeq
It is worth mentioning that the last condition follows since $\la \Gvf, 1 \ra=0$.

Let us define
\beq\notag
\tilde{u}(\xi,\Gt):=(u \circ \Phi)(\xi, \Gt),
\eeq
where $\Phi$ is the conformal mapping defined by \eqnref{def_Phi}.
Since $\Phi$ is conformal, $\tilde{u}$ satisfies
\beq\label{transformed1}
\Delta \tilde{u}(\xi,\Gt) = 0  \quad \mbox{for }(\xi,\Gt)\in\mathbb{R}\times\bigr( (-\pi,-\Gt_0)\cup(-\Gt_0,\Gt_0)\cup(\Gt_0,\pi)\bigr),
\eeq
where $\GD= \p^2/\p\xi^2 + \p^2/\p\Gt^2$. The continuity of the potential in \eqnref{transmission_eigen} is equivalent to
\beq\label{transformed2}
\tilde{u}|_+ =\tilde{u}|_- \quad \mbox{on } \{\Gt=\pm\Gt_0\}.
\eeq
In view of \eqnref{nor_bipolar}, the continuity of the flux is equivalent to
\beq\label{transformed3}
\frac{\p \tilde{u}}{\p \Gt}\Big|_+ - \frac{\p \tilde{u}}{\p \Gt}\Big|_- =\frac{-1}{h(\cdot,\Gt_0)}\tilde{\Gvf}_+ \quad\mbox{on } \{\Gt=\Gt_0\}
\eeq
and
\beq\label{transformed4}
\frac{\p \tilde{u}}{\p \Gt}\Big|_+ - \frac{\p \tilde{u}}{\p \Gt}\Big|_- =\frac{1}{h(\cdot,\Gt_0)}\tilde{\Gvf}_- \quad\mbox{on } \{\Gt=-\Gt_0\},
\eeq
where
\beq\notag
\tilde{\Gvf}_+(\xi)=(\Gvf\circ\Phi)(\xi, \Gt_0) \quad \mbox{and}\quad \tilde{\Gvf}_-(\xi)=(\Gvf\circ\Phi)(\xi, -\Gt_0).
\eeq
Since $|(x_1,x_2)| \to \infty$ if and only if $|(\xi, \Gt)| \to 0$, it follows from the last condition in \eqnref{transmission_eigen} that $\tilde{u}(\xi,\Gt)\to 0$ as $(\xi,\Gt) \to (0,0)$, or simply
\beq\label{transformed5}
\tilde{u}(0,0)= 0.
\eeq
We now impose boundary conditions on $\{\Gt=\pm\pi\}$. Since lines $\{\Gt=\pm\pi\}$ correspond to the line segment $(-\alpha,\alpha)$ in $(x_1,x_2)$-plane, we have
\beq\label{transformed6}
\tilde{u}(\xi,\pi)=\tilde{u}(\xi,-\pi),\quad
\frac{\p \tilde{u}}{\p \Gt}(\xi,\pi)=\frac{\p\tilde{u}}{\p \Gt}(\xi,-\pi)
\quad\mbox{for all } \xi\in \mathbb{R},
\eeq
which implies that $\tilde{u}(\xi, \Gt)$ is periodic in $\Gt$ with the period $2\pi$.

It is straightforward to find the solution $\tilde{u}$ to \eqnref{transformed1} satisfying \eqnref{transformed2}-\eqnref{transformed6}. Denote by $\mathcal{F}$ the Fourier transformation in $\xi$-variable, namely,
\beq\notag
\Fcal[f](s)=\frac{1}{\sqrt{2\pi}}\int_{-\infty}^\infty f(\xi)e^{-i s \xi } d\xi.
\eeq
If we apply $\Fcal$ to the both sides of the equation \eqnref{transformed1} and denote
\beq\notag
v(s, \Gt)= \Fcal[\tilde{u}(\cdot, \Gt)](s), \quad s \in \Rbb,
\eeq
then we have
\beq\notag
\frac{d^2 v}{d \Gt^2} - s^2 v =0.
\eeq
We solve this ordinary differential equation with conditions \eqnref{transformed2}-\eqnref{transformed6}.  In doing so, it is
convenient to consider separately the cases when $\Gvf$ is odd and even with respect to $x_1$-axis.

Suppose that $\Gvf$ is odd with respect to $x_1$-axis, namely
$$
\Gvf(x_1,x_2)=-\Gvf(x_1,-x_2).
$$
Then, one can see easily that
\beq\label{odd}
\tilde{\Gvf}_+= - \tilde{\Gvf}_-,
\eeq
and the solution is given by
\beq\label{tilde_u_sol_odd}
v(s,\Gt) =  \begin{cases}
\ds a(s)  \sinh s(\Gt-\pi), \quad& \Gt_0<\Gt<\pi,\\[0.2cm]
\ds a(s) \frac{\sinh s(\Gt_0-\pi)}{\sinh s\Gt_0}\sinh s\Gt , \quad& -\Gt_0 \le \Gt \le \Gt_0,\\[0.2cm]
\ds a(s) \sinh s(\Gt+\pi), \quad& -\pi<\Gt<-\Gt_0,
\end{cases}
\eeq
where
\beq\notag
a(s)=\frac{\sinh s\Gt_0}{s \sinh s\pi} \Fcal \Big[ \frac{\tilde{\Gvf}_+}{h(\cdot,\Gt_0)} \Big](s).
\eeq

If $\Gvf$ is even with respect to $x_1$-axis, then we have
\beq\label{even}
\tilde{\Gvf}_+= \tilde{\Gvf}_-,
\eeq
and the solution is given by
\beq\label{tilde_u_sol_even}
v(\xi,\Gt) =  (\mbox{const.})+
\begin{cases}
\ds b(s)  \cosh s(\Gt-\pi), \quad& \Gt_0<\Gt<\pi,\\[0.2cm]
\ds b(s) \frac{\cosh s(\Gt_0-\pi)}{\cosh s\Gt_0}\cosh s\Gt , \quad& -\Gt_0 \le \Gt \le \Gt_0,\\[0.2cm]
\ds b(s) \cosh s(\Gt+\pi), \quad& -\pi<\Gt<-\Gt_0,
\end{cases}
\eeq
where
\beq\notag
b(s)=-\frac{\cosh s\Gt_0}{s \sinh s\pi} \Fcal \Big[ \frac{\tilde{\Gvf}_+}{h(\cdot,\Gt_0)} \Big](s).
\eeq

In general case, we decompose $\Gvf$ into the odd and even parts, namely,
\beq\label{oddeven}
\Gvf=\Gvf^o + \Gvf^e
\eeq
where
\beq\notag
\Gvf^o(x_1,x_2):=\frac{1}{2}(\Gvf(x_1,x_2)-\Gvf(x_1,-x_2)) \quad\mbox{and}\quad \Gvf^e(x_1,x_2):=\frac{1}{2}(\Gvf(x_1,x_2)+\Gvf(x_1,-x_2)),
\eeq
and then obtain the solution by superposing solutions corresponding to the odd and even parts, respectively.
By taking the inverse Fourier transform we obtain $\tilde{u}(\xi, \Gt)$, namely,
\beq\label{tildeu}
\tilde{u}(\xi, \Gt)=u(\Phi(\xi, \Gt)) = \Fcal^{-1}[v(\cdot, \Gt)](\xi).
\eeq

We emphasize that
\beq\label{Gvfezero}
\int_{\GG_+} \Gvf^e d\Gs=0.
\eeq
In fact, since $\Gvf \in \Hcal^*$, it holds that
$$
\int_{\GG} \Gvf^o d\Gs + \int_{\GG} \Gvf^e d\Gs=0.
$$
Since $\Gvf^o$ is odd with respect to the $x_1$-axis, we have $\int_{\GG} \Gvf^o d\Gs=0$, and hence $\int_{\GG} \Gvf^e d\Gs=0$. Since $\Gvf^e$ is even, we have \eqnref{Gvfezero}.
%%%%%%%%%%%%%%%%%%%%%%%%%%%%%%%%%%%%%%%%%%%%%%%%%%%%%%%%%%%%%%%%%%
\subsection{Single layer potential on $\p\GO$}
%%%%%%%%%%%%%%%%%%%%%%%%%%%%%%%%%%%%%%%%%%%%%%%%%%%%%%%%%%%%%%%%%%

In particular, we obtain from \eqnref{tildeu} the following formula $\mathcal{S}_{\p\GO}[\Gvf]$ on $\p\GO$:
\beq\label{fourier}
\Scal_{\p\GO}[\Gvf](\Phi(\xi, \pm \Gt_0)) = \Fcal^{-1}[\pm  A^o + A^e](\xi) ,
\eeq
where
\begin{align*}
 A^o(s) &:= \frac{-\sinh s\Gt_0 \sinh s(\pi-\Gt_0) }{s \sinh s\pi} \Fcal\Big[ \frac{\tilde{\Gvf}_+^o}{h(\cdot,\Gt_0)} \Big](s),\\
 A^e(s) &:= \frac{-\cosh s\Gt_0 \cosh s(\pi-\Gt_0)}{s \sinh s\pi}  \Fcal\Big[ \frac{\tilde{\Gvf}_+^e}{h(\cdot,\Gt_0)} \Big](s).
\end{align*}
It is worth reminding here that
\beq\notag
\tilde{\Gvf}_+^o(\xi)= \Gvf^o(\Phi(\xi, \Gt_0)),
\eeq
and likewise for $\tilde{\Gvf}_+^e(\xi)$.

Above formula suggest to introduce an operator $\Ccal$ as follows: for $\Gvf \in \Hcal^*$
\beq\notag
\Ccal[\Gvf](\xi):= \begin{bmatrix} \Ccal_+[\Gvf](\xi) \\ \Ccal_-[\Gvf](\xi) \end{bmatrix}
\eeq
where
\beq\notag
\Ccal_{\pm}[\Gvf](\xi):= \frac{\Gvf(\Phi(\xi, \pm \Gt_0))}{h(\xi,\Gt_0)}.
\eeq
Note that $\Ccal^{-1}$ is given by
\beq\label{Cinverse}
\Ccal^{-1} \begin{bmatrix} g_1 \\ g_2 \end{bmatrix}(x) = \begin{cases} h(\Psi_1(x),\Gt_0) g_1(\Psi_1(x)), \quad& x\in \GG_+, \\
h(\Psi_1(x),\Gt_0) g_2(\Psi_1(x)), \quad & x\in \GG_-. \end{cases}
\eeq
We also define a multiplication operator $\Mcal_h$ by
\beq\notag
\Mcal_h [f]: = h (\xi, \Gt_0)  f(\xi).
\eeq

Let
\begin{align}
p_1(s)&:= \frac{\sinh s\Gt_0 \sinh s(\pi-\Gt_0) }{s \sinh s\pi}, \label{def_pone} \\
p_2(s)&:= \frac{\cosh s\Gt_0 \cosh s(\pi-\Gt_0)}{s \sinh s\pi}. \label{def_ptwo}
\end{align}
Observe that $p_1$ is a smooth function on $\Rbb^1$ and $p_2$ is smooth except at $0$ where the following holds:
\beq\label{ptwozero}
\lim_{s \to \pm 0} |s|^2 p_2(s) = 1.
\eeq
Both of them are positive and satisfy
\beq\label{plimit}
\lim_{s \to \pm \infty} |s| p_j(s) = \frac{1}{2},\quad j=1,2.
\eeq

Then, since $\Gvf_+^o(\xi)= -{\Gvf}_-^o(\xi)$ and ${\Gvf}_+^e(\xi)= {\Gvf}_-^e(\xi)$, we see from \eqnref{fourier} that
\begin{align*}
\Fcal \Mcal_h \Ccal_+ \Scal_{\p\GO}[\Gvf] (s) & = - p_1(s) \Fcal \Ccal_+ [\Gvf^o] (s) - p_2(s) \Fcal\Ccal_+ [\Gvf^e] (s), \\
\Fcal \Mcal_h \Ccal_- \Scal_{\p\GO}[\Gvf] (s) & = - p_1(s) \Fcal\Ccal_- [\Gvf^o] (s) - p_2(s) \Fcal\Ccal_- [\Gvf^e] (s).
\end{align*}
Let
\beq\label{deff}
\Bf(s) = \begin{bmatrix} f_1(s) \\ f_2(s) \end{bmatrix} = \begin{bmatrix} \Fcal\Ccal_+ [\Gvf^o] (s) \\[1mm]
\Fcal\Ccal_+ [\Gvf^e] (s) \end{bmatrix}.
\eeq
Then we have
$$
f_1 + f_2 = \Fcal\Ccal_+ [\Gvf] \quad\mbox{and}\quad -f_1 + f_2 = \Fcal\Ccal_- [\Gvf],
$$
in other words,
\beq\label{Gvfform}
\Gvf = 2 \Ccal^{-1} \Fcal^{-1} \GL^{-1} \begin{bmatrix} f_1 \\ f_2 \end{bmatrix},
\eeq
where
\beq\label{jay}
\GL= \begin{bmatrix} 1 & -1 \\ 1 & 1 \end{bmatrix}.
\eeq
So, we obtain
\begin{align*}
2\Fcal \Mcal_h \Ccal \Scal_{\p\GO} (\GL\Fcal\Ccal)^{-1} \begin{bmatrix} f_1 \\ f_2 \end{bmatrix}
& = \begin{bmatrix} - p_1 f_1 -  p_2 f_2 \\ p_1 f_1 - p_2 f_2 \end{bmatrix} .
\end{align*}
Note that
$$
\begin{bmatrix} - p_1 f_1 -  p_2 f_2 \\ p_1 f_1 - p_2 f_2 \end{bmatrix} = - 2 \GL^{-1} P(s) \Bf(s)
$$
where
\beq\notag
P= \begin{bmatrix} p_1 & 0 \\ 0 & p_2 \end{bmatrix} .
\eeq

If we define the transformed single layer potential $\Sbb$ by
\beq\label{Sbb}
\Sbb := \GL \Fcal \Mcal_h \Ccal \Scal_{\p\GO} (\GL\Fcal\Ccal)^{-1},
\eeq
it follows that
\beq\label{SbbBf}
\Sbb [\Bf](s) = - P(s) \Bf(s).
\eeq

%%%%%%%%%%%%%%%%%%%%%%%%%%%%%%%%%%%%%%%%%%%%%%%%%%%%%%%%%%%%%%%%%%
\subsection{Inner product}
%%%%%%%%%%%%%%%%%%%%%%%%%%%%%%%%%%%%%%%%%%%%%%%%%%%%%%%%%%%%%%%%%%

Observe that if we make changes of variables $x=\Phi(\xi, \pm \Gt_0)$, then we see from \eqnref{jacobian} that
$$
\int_{\p\GO} \Gvf(x) \overline{\psi(x)} d\Gs(x) = \int_{\Rbb} \Ccal_+[\Gvf](\xi) \Mcal_h \overline{\Ccal_+[\psi](\xi)} d\xi
+ \int_{\Rbb} \Ccal_-[\Gvf](\xi) \Mcal_h \overline{\Ccal_-[\psi](\xi)} d\xi.
$$
It then follows from Parseval's identity that
\begin{align*}
\int_{\p\GO} \Gvf(x) \overline{\psi(x)} d\Gs(x)
& = \int_{\Rbb} \Fcal \Ccal_+[\Gvf](s) \overline{\Fcal \Mcal_h \Ccal_+[\psi](s)} ds + \int_{\Rbb} \Fcal \Ccal_-[\Gvf](s) \overline{\Fcal\Mcal_h \Ccal_-[\psi](s)} ds \\
& = \int_{\Rbb} \Fcal \Ccal[\Gvf](s) \cdot \overline{\Fcal \Mcal_h \Ccal[\psi](s)} ds .
\end{align*}
Since $\GL^{-1}= \frac{1}{2} \GL^T$, we have
\begin{align*}
\int_{\p\GO} \Gvf(x) \overline{\psi(x)} d\Gs(x)
= \frac{1}{2} \int_{\Rbb} \GL \Fcal \Ccal[\Gvf](s) \cdot \overline{\GL \Fcal \Mcal_h \Ccal[\psi](s)} ds .
\end{align*}

Let
\beq
\Bf=(f_1,f_2)^T:= \GL \Fcal \Ccal[\Gvf] \quad\mbox{and}\quad \Bg=(g_1,g_2)^T:= \GL \Fcal \Ccal[\psi].
\eeq
Then, we have from \eqnref{SbbBf}
\begin{align*}
\la \Gvf, \psi\ra_{\Hcal^*} &= - \la \Gvf, \Scal_{\p\GO}[\psi]\ra \\
&= - \frac{1}{2} \int_{\Rbb} \GL \Fcal \Ccal[\Gvf](s) \cdot \overline{\GL \Fcal \Mcal_h \Ccal \Scal_{\p\GO}[\psi](s)} ds \\
& = - \frac{1}{2} \int_{\Rbb} \Bf(s) \cdot \overline{\Sbb[\Bg](s)} ds \\
& = \frac{1}{2} \int_{\Rbb} \Bf(s) \cdot P(s) \overline{\Bg(s)} ds.
\end{align*}

Let
\beq\label{newinner}
\la \Bf, \Bg \ra_{W^*} : = \frac{1}{2} \int_{\Rbb} \Bf(s) \cdot P(s) \overline{\Bg(s)}\; ds = \frac{1}{2} \int_{\Rbb} \left[ p_1(s) f_1(s)  \overline{g_1(s)} + p_2(s) f_2(s)  \overline{g_2(s)} \right] ds,
\eeq
and let $\| \cdot \|_{W^*}$ be the norm induced by this inner product. Let $W^*$ be the space of $\Bf$ such that $\| \Bf \|_{W^*} < \infty$. In fact, $W^*= L^2(\Rbb, p_1) \times L^2(\Rbb, p_2)$, where $L^2(\Rbb, p_j)$ is the $L^2$ space weighted by $p_j$. It is worth mentioning that because of \eqnref{plimit} $L^2(\Rbb, p_1)$ is the Fourier transform of the space $H^{-1/2}(\Rbb)$. On the other hand, because of \eqnref{ptwozero}, if $g \in L^2(\Rbb, p_2)$, $g(0)=0$ if g is continuous at $0$. So, $L^2(\Rbb, p_2)$ is the Fourier transform of $H^{-1/2}_0(\Rbb)$. This is natural in view of \eqnref{Gvfezero}.

Define
\beq\label{Udef}
U:= \GL\Fcal\Ccal.
\eeq
We have shown that $U$ is a unitary transformation from $\Hcal^*$ onto $W^*$, i.e.,
\beq\label{unitary}
\la \Gvf, \psi \ra_{\Hcal^*} = \la U[\Gvf], U[\psi] \ra_{W^*}.
\eeq

%%%%%%%%%%%%%%%%%%%%%%%%%%%%%%%%%%%%%%%%%%%%%%%%%%%%%%%%%%%%%%%%%%
\subsection{Neumann-Poincar\'e operator}
%%%%%%%%%%%%%%%%%%%%%%%%%%%%%%%%%%%%%%%%%%%%%%%%%%%%%%%%%%%%%%%%%%

We now compute the NP operator and show that
\beq\label{NPform}
\Kcal^*_{\p\GO}[\Gvf](\Phi(\xi, \pm \Gt_0)) = h(\xi,\Gt_0)\, \Fcal^{-1}[\pm B^o - B^e](\xi) ,
\eeq
where
\begin{align*}
B^o(s)&:= \eta(s) \Fcal \Big[ \frac{\tilde{\Gvf}_+^o}{h(\cdot,\Gt_0)} \Big](s),\\
B^e(s)&:= \eta(s) \Fcal \Big[ \frac{\tilde{\Gvf}_+^e}{h(\cdot,\Gt_0)} \Big](s).
\end{align*}
Here, $\eta$ is defined by
\beq\label{etadef}
\eta(s):=\frac{1}{2}\frac{\sinh s(\pi-2\Gt_0)}{\sinh s\pi}.
\eeq

We only show \eqnref{NPform} for $\Gt=\Gt_0$. Let $u=\Scal_{\p\GO}[\Gvf]$ and $\tilde{u}(\xi, \Gt)=\Scal_{\p\GO}[\Gvf](\Phi(\xi, \Gt))$ as before. We see from \eqnref{singlejump} and \eqnref{nor_bipolar} that
\begin{align}
\Kcal^*_{\p\GO}[\Gvf]\big( \Phi (\xi,  \Gt_0)\big)
&=\frac{1}{2}\left( \frac{\p u}{\p\nu}\Big|_+ + \frac{\p u}{\p\nu}\Big|_- \right)\big( \Phi (\xi,  \Gt_0)\big) \nonumber \\
&= - \frac{h(\xi,\Gt_0)}{2}\left(\frac{\p \tilde{u}}{\p\Gt}\Big|_+ (\xi,  \Gt_0) +
\frac{\p \tilde{u}}{\p\Gt}\Big|_- (\xi,  \Gt_0) \right). \label{K_du}
\end{align}
Recall that the region $\mathbb{R}^2\setminus \overline{\Omega}$ outside the intersecting disks corresponds to the region $\{|\theta|<\theta_0\}$. So the subscript $+$ in $\frac{\p}{\p\theta}|_\pm$ indicates the limit from the region $\{|\theta|<\theta_0\}$. Similarly, the subscript $-$ indicates the limit from the region $\{\theta_0<|\theta|<\pi\}$.

Let $v=v^o+v^e:= \Fcal[\tilde{u}]$. One can easily see from \eqnref{tilde_u_sol_odd} and \eqnref{tilde_u_sol_even} that
\begin{align*}
\frac{\p v^o}{\p\Gt}\Big|_+(s, \Gt_0)
&=\frac{\sinh s(\Gt_0-\pi)\cosh s\Gt_0}{\sinh s\pi}\hskip .5mm \Fcal\Big[ \frac{\tilde{\Gvf}_+^o}{h(\cdot,\Gt_0)}\Big](s),  \\
\frac{\p v^o}{\p\Gt}\Big|_- (s,\Gt_0)
&=\frac{\cosh s(\Gt_0-\pi)\sinh s\Gt_0}{\sinh s\pi} \hskip .5mm\Fcal\Big[ \frac{\tilde{\Gvf}_+^o}{h(\cdot,\Gt_0)}\Big](s), \\
\frac{\p v^e}{\p\Gt}\Big|_+ (s,\Gt_0)
&= - \frac{\cosh s(\Gt_0-\pi)\sinh s\Gt_0}{\sinh s\pi} \hskip .5mm \Fcal\Big[ \frac{\tilde{\Gvf}_+^e}{h(\cdot,\Gt_0)}\Big](s), \\
\frac{\p v^e}{\p\Gt}\Big|_- (s,\Gt_0)
&=-\frac{\sinh s(\Gt_0-\pi)\cosh s\Gt_0}{\sinh s\pi}\hskip .5mm \Fcal\Big[ \frac{\tilde{\Gvf}_+^e}{h(\cdot,\Gt_0)}\Big](s).
\end{align*}
Since
$$
\sinh s(\Gt_0-\pi)\cosh s\Gt_0+\cosh s(\Gt_0-\pi)\sinh s\Gt_0=-\sinh s(\pi-2\Gt_0),
$$
we obtain
$$
\frac{\p v}{\p\Gt}\Big|_+ (s,\Gt_0) + \frac{\p v}{\p\Gt}\Big|_- (s,\Gt_0)
=-\frac{\sinh s(\pi-2\Gt_0)}{\sinh s\pi}\left(\Fcal\Big[ \frac{\tilde{\Gvf}_+^o}{h(\cdot,\Gt_0)}\Big](s) -
\Fcal\Big[ \frac{\tilde{\Gvf}_+^e}{h(\cdot,\Gt_0)}\Big](s)\right).
$$
So, we obtain \eqnref{NPform} for $\Gt=\Gt_0$.

Define $\Bf$ by \eqnref{deff}. Then, \eqnref{NPform} can be written as
\beq\label{NPform2}
2\Fcal\Ccal \Kcal^*_{\p\GO} (\GL\Fcal\Ccal)^{-1} [\Bf] = \eta \begin{bmatrix} f_1- f_2 \\ - f_1-f_2 \end{bmatrix} = 2\eta \GL^{-1} \begin{bmatrix} f_1 \\ -f_2 \end{bmatrix}.
\eeq
Let
\beq\label{Kbbdef}
\Kbb^*:= U \Kcal_{\p\GO}^* U^{-1} = \GL \Fcal\Ccal \Kcal_{\p\GO}^* (\GL \Fcal\Ccal)^{-1}.
\eeq
Then using  \eqnref{Gvfform}, \eqnref{NPform} can be rewritten as
\beq\label{Kbbeta}
\Kbb^* [\Bf](s) = \eta(s) \begin{bmatrix} f_1(s) \\ -f_2(s) \end{bmatrix}
\eeq
on $W^*$.

Note that $\eta$ is monotonically decreasing for $s \ge 0$, tends to $0$ as $s \to \infty$, and satisfies
\beq\label{etazero}
\eta(s) = b \left[ 1 - \frac{2\Gt_0(\pi-\Gt_0)}{3} s^2 + O(s^4) \right]
\eeq
where
\beq\label{bpo}
b := \eta(0)= \frac{1}{2} - \frac{\Gt_0}{\pi} .
\eeq
We can infer from \eqnref{Kbbeta} that $[-b, b]$ is the spectrum of $\Kbb^*$ (and hence of $\Kcal_{\p\GO}^*$),  and it is a continuous spectrum. We will investigate the nature of the spectrum in the following section.

%%%%%%%%%%%%%%%%%%%%%%%%%%%%%%%%%%%%%%%%%%%%%%%%%%%%%%%%%%%%%%%%%%%%%%%%%%%%%%%%
\section{Spectral resolution of NP-operator}\label{sec4}
%%%%%%%%%%%%%%%%%%%%%%%%%%%%%%%%%%%%%%%%%%%%%%%%%%%%%%%%%%%%%%%%%%%%%%%%%%%%%%%%

%%%%%%%%%%%%%%%%%%%%%%%%%%%%%%%%%%%%%%%%%%%%%%%%%%%
\subsection{Resolution of the identity}
%%%%%%%%%%%%%%%%%%%%%%%%%%%%%%%%%%%%%%%%%%%%%%%%%%%

Let $\chi_s$ be the characteristic function on $(-\infty, s]$ for each $s\in \Rbb$. Let
$$
\Pi_1 \Bf = \begin{bmatrix} f_1 \\ 0 \end{bmatrix} \quad\mbox{and}\quad \Pi_2 \Bf = \begin{bmatrix} 0 \\ f_2 \end{bmatrix}.
$$
We define a pair of operators $\Pcal_s^1$ and $\Pcal_s^2$ on $W^*$ by
\beq\notag
\Pcal_s^j [\Bf] = \chi_s \Pi_j \Bf, \quad j=1,2.
\eeq
One can see that $\Pcal_s^j$ is an orthogonal projection, namely, $\Pcal_s^j$ is self-adjoint and $(\Pcal_s^j)^2=\Pcal_s^j$. We also have
\beq\label{proj2}
\Pcal_s^1 \Pcal_t^2 = \Pcal_t^2 \Pcal_s^1 = 0 \quad\mbox{for all } s, t \in \Rbb.
\eeq
One can also see that the following holds:
\beq\label{proj3}
\begin{cases}
\ds \lim_{s \to \infty} \Pcal_s^j [\Bf] = \Pi_j \Bf ,  \\
\ds \lim_{s \to -\infty} \Pcal_s^j [\Bf] =  0, \\[0.2cm]
\ds \lim_{t \to s} \Pcal_t^j [\Bf] = \Pcal_s^j [\Bf],  \\[0.2cm]
\ds \Pcal_t^j \Pcal_s^j [\Bf] = \Pcal_{t \wedge s}^j [\Bf] ,
\end{cases}
\eeq
where $t\wedge s$ denotes the minimum of $t$ and $s$, and the convergence is in $W^*$. As a consequence, we have
\beq\label{proj4}
\lim_{s \to \infty} (\Pcal_s^1 [\Bf] + \Pcal_s^2 [\Bf]) = \Bf.
\eeq

We see from \eqnref{newinner} that
\beq\label{proj5}
\la \Bf, \Pcal_s^j[\Bf] \ra_{W^*} = \frac{1}{2} \int_{-\infty}^s p_j |f_j|^2 ds ,\quad j=1,2.
\eeq
It then follows that
\begin{align}
\la \Bf, \Bf \ra_{W^*} &= \frac{1}{2} \int_{\Rbb} \left[ p_1 |f_1|^2 + p_2 |f_2|^2 \right] \, ds. \nonumber \\
&= \int_{-\infty}^\infty \frac{d}{ds} \big\la \Bf, \Pcal_s^1 [\Bf] \big\ra_{W^*} \,ds + \int_{-\infty}^\infty
\frac{d}{ds} \big\la \Bf, \Pcal_s^2 [\Bf] \big\ra_{W^*} \,ds. \label{resolid}
\end{align}
By the change of variables $t=\eta(s)$ for $s \ge 0$ and $t=\eta(-s)$ for $s<0$, we obtain
$$
\int_{-\infty}^\infty \frac{d}{ds} \big\la \Bf, \Pcal_s^1 [\Bf] \big\ra_{W^*} \,ds
= \int_{0}^{b}
d \left\la \Bf,  \left( \Pcal_{-\eta^{-1}(t)}^1 - \Pcal_{\eta^{-1}(t)}^1 \right)
[\Bf] \right\ra_{W^*} ,
$$
which can be written as
$$
\int_{-\infty}^\infty \frac{d}{ds} \big\la \Bf, \Pcal_s^1 [\Bf] \big\ra_{W^*} \,ds
= \int_{0}^{b}
d \left\la \Bf,  \left( \Pcal_{-\eta^{-1}(t)}^1 - \Pcal_{\eta^{-1}(t)}^1 + \Ical\right)
[\Bf] \right\ra_{W^*} ,
$$
where $\Ical$ is the identity operator on $W^*$. If we make change of variables $-t=\eta(s)$ for $s \ge 0$ and $-t=\eta(-s)$ for $s<0$, then we have
$$
\int_{-\infty}^\infty
\frac{d}{ds} \big\la \Bf, \Pcal_s^2 [\Bf] \big\ra_{W^*} \,ds = \int_{-b}^{0}
d \left\la \Bf,  \left(  \Pcal_{\eta^{-1}(t)}^2 - \Pcal_{-\eta^{-1}(t)}^2\right)
[\Bf] \right\ra_{W^*} .
$$
So we obtain from \eqnref{resolid} that
\begin{align}
\la \Bf,\Bf\ra_{W^*}
&= \int_{0}^{b}
d \left\la \Bf,  \left( \Pcal_{-\eta^{-1}(t)}^1 - \Pcal_{\eta^{-1}(t)}^1 + \Ical\right)
[\Bf] \right\ra_{W^*} \nonumber \\
&\quad\qquad + \int_{-b}^{0}
d \left\la \Bf, \left( \Pcal_{\eta^{-1}(t)}^2 - \Pcal_{-\eta^{-1}(t)}^2\right)
[\Bf] \right\ra_{W^*} . \label{resolid2}
\end{align}

We now define a family of projection operators $\{\Ebb_t\}_{t\in[-b,b]}$ on $W^*$ as follows:
\beq\label{etdef}
\Ebb_t :=
\begin{cases}
 \Pcal^1_{-\eta^{-1}(t)}-\Pcal^1_{\eta^{-1}(t)}+\Ical \, , \quad& t\in (0,b],\\[0.5em]
 \Pcal_{\eta^{-1}(-t)}^2 - \Pcal_{-\eta^{-1}(-t)}^2 \, , \quad& t\in [-b,0),
\end{cases}
\eeq
and
\beq\label{etzero}
\Ebb_0 := \lim_{t \to 0+} \Ebb_t,
\eeq
in other words,
\beq\notag
\Ebb_0[\Bf] = \lim_{t \to 0+} \Ebb_t [\Bf]
\eeq
for $\Bf \in W^*$.
In view of the first two identities in \eqnref{proj3}, we have
\beq\label{etzero2}
\Ebb_0[\Bf] = - \Pi_1 \Bf +\Bf = \Pi_2 \Bf.
\eeq
Now one can easily see from \eqnref{proj2}, \eqnref{proj3}, and \eqnref{etzero} that
the following holds:
\beq\label{etes}
\begin{cases}
\ds\Ebb_t\Ebb_s = \Ebb_{t \wedge s}, \\
\ds \Ebb_{b} [\Bf] = \Bf, \quad \Ebb_{-b} [\Bf] = 0, \\
\ds \lim_{t \to s+} \Ebb_t [\Bf] = \Ebb_s [\Bf]
\end{cases}
\eeq
for all $s,t \in [-b, b]$.

We obtain the following proposition.

\begin{prop}\label{resol_prop}
The family $\{\Ebb_t\}_{t\in[-b,b]}$ is a resolution of the identity on $W^*$, namely, it holds that
\beq \label{resol_iden_eq1}
\la \Bg,\Bf \ra_{W^*}=\int_{-b}^{b} d \la \Bg,\Ebb_t[\Bf] \ra_{W^*}
\eeq
for all $\Bg,\Bf\in W^*$. Moreover, we have
\beq\label{etconti}
\lim_{t \to s} \Ebb_t [\Bf] = \Ebb_s [\Bf]
\eeq
for all $s \in [-b, b]$.
\end{prop}
\pf
That $\{\Ebb_t\}_{t\in[-b,b]}$ is a resolution of the identity on $W^*$ is an immediate consequence of \eqnref{etes}. To prove \eqnref{etconti}, it suffices to consider the case when $s=0$ because of the third identity in \eqnref{proj3}. But one can see from the definition \eqnref{etdef} that
$$
\lim_{t \to 0-} \Ebb_t [\Bf] =
\lim_{t \to \infty} \Pcal_{t}^2[\Bf] - \lim_{t \to -\infty} \Pcal_{t}^2 [\Bf] = \Pi_2 \Bf .
$$
So by \eqnref{etzero2} we have \eqnref{etconti} for $s=0$, and the proof is complete.
\qed

%%%%%%%%%%%%%%%%%%%%%%%%%%%%%%%%%%%%%%%%%%%%%%%%%%%%%%%%%%%%%%%%%%%%%
\subsection{Spectral resolution of the NP-operator}
%%%%%%%%%%%%%%%%%%%%%%%%%%%%%%%%%%%%%%%%%%%%%%%%%%%%%%%%%%%%%%%%%%%%%

We now derive the spectral resolution of $\Kbb^*$.

\begin{theorem}\label{main_thm}
Let $\{\Ebb_t\}$ be the resolution of the identity defined in {\rm\eqnref{etdef}}. Then
we have the following spectral resolution of $\Kbb^* $ on $W^*$:
\beq\label{specresol2}
\Kbb^* = \int_{-b}^{b}
t \, d  \Ebb_t.
\eeq
In other words, it holds that
\beq\label{specresol}
\la \Bg,\Kbb^* [\Bf]\ra_{W^*} = \int_{-b}^{b}
t \, d \la \Bg , \Ebb_t[\Bf]\ra_{W^*}
\eeq
for all $\Bg,\Bf \in W^*$.
\end{theorem}

\pf
It suffices to show \eqnref{specresol} when $\Bg=\Bf$ by the polarization identity. We obtain from \eqnref{newinner} and \eqnref{Kbbeta}
$$
\la \Bf,\Kbb^* [\Bf]\ra_{W^*} =
\frac{1}{2} \int_{\Rbb} \eta \left[ p_1 |f_1|^2 - p_2 |f_2|^2 \right] \, ds.
$$
It then follows from \eqnref{proj5} and \eqnref{etadef} that
\begin{align*}
\la \Bf,\Kbb^* [\Bf]\ra_{W^*}
= \int_{-\infty}^\infty \eta(s) \frac{d}{ds} \big\la \Bf, \Pcal_s^1[\Bf] \big\ra_{W^*} \,ds
- \int_{-\infty}^\infty \eta(s) \frac{d}{ds} \big\la \Bf, \Pcal_s^2 [\Bf] \big\ra_{W^*} \,ds.
\end{align*}
By applying the same changes of variables as before, $t=\eta(|s|)$ or $t=-\eta(|s|)$, we obtain
\begin{align*}
\la \Bf,\Kbb^* [\Bf]\ra_{W^*}
&= \int_{0}^{b} t\,
d \left\la \Bf, \left(\Pcal_{-\eta^{-1}(t)}^+ - \Pcal_{\eta^{-1}(t)}^+ + \Ical\right)
\Bf \right\ra_{W^*}\\
&\qquad +
\int_{-b}^{0} t\,
d \left\la \Bf, \left( \Pcal_{\eta^{-1}(-t)}^- - \Pcal_{-\eta^{-1}(-t)}^-\right)
\Bf \right\ra_{W^*} \\
&= \int_{0}^{b} t\, d \la \Bf, \Ebb_t \Bf \ra_{W^*} +
\int_{-b}^{0} t\, d \la \Bf, \Ebb_t \Bf \ra_{W^*}.
\end{align*}
Since $\Ebb_t$ is continuous at $0$, we obtain \eqnref{specresol} when $\Bg=\Bf$. This completes the proof.
\qed

%%%%%%%%%%%%%%%%%%%%%%%%%%%%%%%%%%%%%%%%%%%%%%%%%%%%%%%%%%%%%%%%%%%%%%
\subsection{Spectrum of the NP operator on the intersecting disks}
%%%%%%%%%%%%%%%%%%%%%%%%%%%%%%%%%%%%%%%%%%%%%%%%%%%%%%%%%%%%%%%%%%%%%%

Recall that the point spectrum $s$ is characterized by the condition (see \cite{Yo})
$$
\lim_{t \to s-} \Ebb_t \neq \Ebb_s.
$$
Since $\Ebb_t$ is continuous for all $s$, namely, \eqnref{etconti} holds, we infer that there is only continuous spectrum and no point spectrum of $\Kbb^*$. In this section we prove that there is only absolutely continuous spectrum. For that we derive the following lemma.

\begin{lemma}\label{specmeasure}
For all $\Bf, \Bg \in W^*$, it holds that
\beq\label{specder}
\frac{d}{dt} \la \Bf , \Ebb_t[\Bg]\ra_{W^*} =
\begin{cases}
\ds \frac{p_1(s)}{2|\eta'(s)|} \Bigr[f_1(s) \overline{g_1(s)} + f_1(-s) \overline{g_1(-s)}\Bigr], \quad& t\in (0,b],  \\[0.3cm]
\ds \frac{p_2(s)}{2|\eta'(s)|} \Bigr[f_2(s) \overline{g_2(s)} + f_2(-s) \overline{g_2(-s)}\Bigr], \quad& t\in [-b, 0),
\end{cases}
\eeq
where $s=\eta^{-1}(t)$ if $t >0$, $s=\eta^{-1}(-t)$ if $t <0$ and
the equality holds in $L^1([-b, b])$.
\end{lemma}
\pf
By the polarization identity it suffices to consider the case when $\Bf=\Bg$.
By the definition \eqnref{etdef} of the resolution of identity $\{\Ebb_t\}$, we have
$$
\la \Bf , \Ebb_t[\Bf]\ra_{W^*} :=
\begin{cases}
 \la\Bf, \Pcal^1_{-\eta^{-1}(t)}[\Bf]\ra_{W^*}
-\la\Bf,\Pcal^1_{\eta^{-1}(t)}[\Bf]\ra_{W^*}+
\la\Bf,\Bf\ra_{W^*} \, , \quad& t\in (0,b],\\[0.5em]
 \la\Bf, \Pcal^2_{\eta^{-1}(-t)}[\Bf]\ra_{W^*}
-\la\Bf,\Pcal^2_{-\eta^{-1}(-t)}[\Bf]\ra_{W^*}
 \, , \quad& t\in [-b,0).
\end{cases}
$$
So, if $t>0$, we have from \eqnref{proj5} that
$$
\la \Bf , \Ebb_t[\Bf]\ra_{W^*} =\big\la \Bf, \Bf \big\ra_{W^*}- \frac{1}{2} \int_{-\eta^{-1}(t)}^{\eta^{-1}(t)}
p_1 |f_1|^2 \,ds.
$$
Since $p_1$ is even, we obtain
$$
\frac{d}{dt} \la \Bf , \Ebb_t[\Bf]\ra_{W^*} = -\frac{p_1(\eta^{-1}(t))}{2\eta'(\eta^{-1}(t))} \Bigr[ |f_1 (\eta^{-1}(t))|^2 + |f_1 (-\eta^{-1}(t))|^2 \Bigr] .
$$
Since $\eta'$ is negative, we obtain \eqnref{specder} for $t>0$. The case when $t<0$ can be proved similarly. \qed

One can easily see that $\frac{d}{dt} \la \Bf , \Ebb_t[\Bf]\ra_{W^*}$ is integrable on $[-b, b]$. In fact, since $\Bf \in W^*$, we have
\begin{align*}
\int_{-b}^{b} \left| \frac{d}{dt} \la \Bf , \Ebb_t[\Bf]\ra_{W^*} \right| dt
&= \frac{1}{2} \int_0^{\infty} p_1(s) \left( |f_1(s)|^2 + |f_1(-s)|^2 \right) ds \\
& \qquad + \frac{1}{2} \int_0^{\infty}  p_2(s) \left( |f_2(s)|^2 + |f_2(-s)|^2 \right) ds < \infty .
\end{align*}
Since $\Kbb^*$ is unitarily equivalent to $\Kcal_{\p\GO}^*$, we obtain the following theorem.

\begin{theorem}\label{spec}
Let $\Gs_{ac}(\Kcal^*_{\p\GO})$, $\Gs_{sc}(\Kcal^*_{\p\GO})$ and $\Gs_{sc}(\Kcal^*_{\p\GO})$ be the absolutely continuous, singularly continuous, and pure point spectrum of $\Kcal^*_{\p\GO}$ on $\Hcal^*$, respectively. Then we have
\beq\notag
\Gs_{ac}(\Kcal^*_{\p\GO})= [-b,b] ,  \quad \Gs_{sc}(\Kcal^*_{\p\GO})= \emptyset, \quad \Gs_{pp}(\Kcal^*_{\p\GO})= \emptyset .
\eeq
\end{theorem}

\begin{figure}[!ht]
\begin{center}
\epsfig{figure=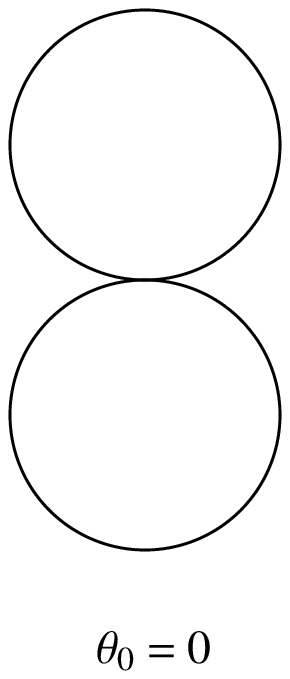,width=4.8cm}
\epsfig{figure=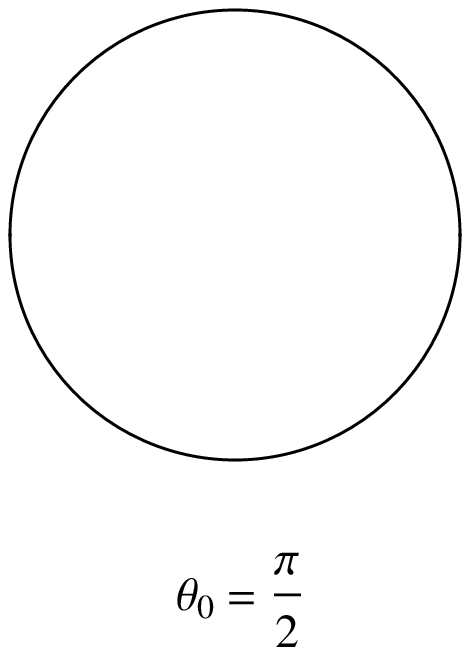,width=4.8cm}
\epsfig{figure=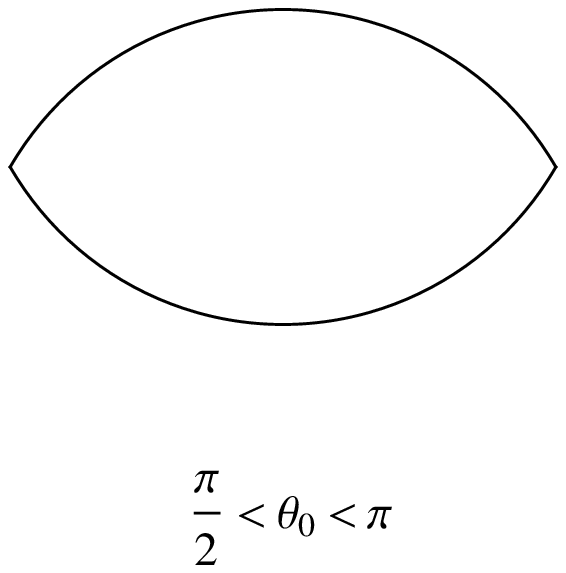,width=4.8cm}
\end{center}
\caption{Disks of various intersecting angles}\label{touching}
\end{figure}

A few remarks on Theorem \ref{spec} are in order. If $\Gt_0=\pi/2$ or $\GO$ is a disk, then $b=0$ which is in accordance with the fact that $\Kcal^*_{\p\GO}=0$ on $\Hcal^*$ when $\GO$ is a disk. It is also interesting to observe what happens when $\Gt_0 \to 0$, or when the intersecting disks approaches to the touching disks. See Fig.\;\ref{touching}. In this case the spectral bound tends to $1/2$. This fact was also observed in the other direction in \cite{BT, BT2, LY1}. It is proved that if two separating disks approach to each other, then there are eigenvalues approaching to $1/2$. It should be mentioned that the NP operator for the touching disk is a hyper-singular operator due to the cusp, and is not well-defined on $H^{-1/2}$.

The analysis of this paper goes through with slight modification for the case when $\Gt_0 > \pi/2$ (see Fig.\;\ref{touching}). In this case the spectral bounds is given by
\beq\label{bounds}
b = \frac{1}{2} - \frac{\pi - \Gt_0}{\pi}= \frac{\Gt_0}{\pi} - \frac{1}{2}.
\eeq
Recently a bound on the essential spectrum of the NP operator on curvilinear polygonal domains has been obtained in \cite{PP}. It is quite interesting to observe that the bound in \eqnref{bounds} is exactly the one obtained there for the intersecting disks.

%%%%%%%%%%%%%%%%%%%%%%%%%%%%%%%%%%%%%%%%%%%%%%%%%%%%%%%%%%%%%%%%%%%%
\section{Plasmonic resonance on the intersecting disks}
%%%%%%%%%%%%%%%%%%%%%%%%%%%%%%%%%%%%%%%%%%%%%%%%%%%%%%%%%%%%%%%%%%%

In this section we consider plasmonic resonance on the intersecting disks.
Suppose that the intersecting disks has dielectric constant $\Ge_c+i \Gd$ ($\Gd$ is the dissipation) and the matrix $\Rbb^2 \setminus \overline{\GO}$ has the dielectric constant  $\Ge_m$. So, the distribution of the dielectric constant is given by
$$
\Ge = (\Ge_c+i\Gd) \chi(\GO) + \Ge_m \chi(\Rbb^2 \setminus \overline{\GO}),
$$
where $\chi$ denotes the indicator function. After normalization we assume that $\Ge_m=1$. To investigate the resonance on the plasmonic structure, we consider the following problem:
\beq\label{trans}
\left\{
\begin{array}{ll}
\nabla \cdot \Ge \nabla u = f \quad &\mbox{in } \Rbb^2, \\
u(x) = O(|x|^{-1}) \quad &\mbox{as } |x| \to \infty,
\end{array}
\right.
\eeq
where $f$ is a source function compactly supported in $\Rbb^2 \setminus \overline{\GO}$ and satisfying
\beq\label{intzero}
\int_{\Rbb^d} f =0.
\eeq
A typical such source functions are polarized dipoles, namely, $f(x)= a \cdot \nabla \Gd_z(x)$ for some $z \in \Rbb^2 \setminus \overline{\GO}$, where $a$ is a constant vector and $\Gd_z$ is the Dirac mass. It is worth mentioning that \eqnref{intzero} is necessary for a solution to \eqnref{trans} to exist. For a given $\Ge_c$ resonance is characterized by the fact
\beq\label{reoinf1}
\| \nabla u_\Gd \|_{L^2(\GO)} \to \infty \quad\mbox{as } \Gd \to 0.
\eeq

Let
\beq\label{lambda}
\Gl:= \frac{\Ge_c +1 + i\Gd}{2(\Ge_c-1) + 2i\Gd}.
\eeq
Then, the solution $u_\Gd$ to \eqnref{trans} can be represented as
\beq\label{solrep}
u_\Gd(x) = q(x) + \Scal_{\p\GO} [\Gvf_\Gd](x), \quad x \in \Rbb^2.
\eeq
where $q$ is the Newtonian potential of $f$, {\it i.e.},
\beq\notag
q(x)= \int_{\Rbb^d} \GG(x-y) f(y) dy,
\eeq
and the potential $\Gvf_\Gd \in \Hcal^*(\p\GO)$ is the solution to the integral equation
\beq\label{solint}
( \Gl I - \Kcal_{\p\GO}^* )[\Gvf_\Gd] = \p_\nu q \quad\mbox{on } \p\GO.
\eeq
See \cite{KS2000} for the proof of this representation.

Let
\beq
\Ecal_t := U^{-1} \Ebb_t U
\eeq
where $\Ebb_t$ is the resolution identity on $W^*$ (Proposition \ref{resol_prop}) and $U$ is the unitary operator in \eqnref{Udef}.
Then we see from Theorem \ref{main_thm} that
\beq
\Kcal^*_{\p\GO} = \int_{-b}^{b}
t \, d  \Ecal_t.
\eeq
We then obtain from \eqnref{solint}
\beq\label{specGvf}
\Gvf_\Gd = \int_{-b}^{b} \frac{1}{\Gl-t} \, d \Ecal_t[\p_\nu q].
\eeq
Since $\| \nabla q \|_{L^2(\GO)} < \infty$, we infer from \eqnref{solrep} that \eqnref{reoinf1} is equivalent to
\beq\label{reoinf2}
\| \nabla \Scal_{\p\GO}[\Gvf_\Gd] \|_{L^2(\GO)} = \| \nabla (u_\Gd-q) \|_{L^2(\GO)} \to \infty \quad\mbox{as } \Gd \to 0.
\eeq
We see from \eqnref{singlejump} that
\begin{align*}
\| \nabla \Scal_{\p\GO}[\Gvf_\Gd] \|_{L^2(\GO)}^2
&= \int_{\p\GO} \pd{}{\nu} \Scal_{\p\GO}[\Gvf_\Gd] \Big|_{-} \overline{\Scal_{\p\GO}[\Gvf_\Gd]} \, d\Gs \\
&= \big\la \Gvf_\Gd , \big(\frac{1}{2}I - \Kcal_{\p\GO}^*\big)[\Gvf_\Gd] \big\ra_{\Hcal^*} .
\end{align*}
Using \eqnref{specGvf} we see that
$$
\big\la \Gvf_\Gd , \big(\frac{1}{2}I - \Kcal_{\p\GO}^*\big)[\Gvf_\Gd] \big\ra_{\Hcal^*}
= \int_{-b}^{b} \left( \frac{1}{2}-t \right) \, d \la \Gvf_\Gd, \Ecal_t[\Gvf_\Gd] \ra_{\Hcal^*}.
$$
So, we infer that
$$
C_1 \| \Gvf_\Gd \|_{\Hcal^*} \le \| \nabla \Scal_{\p\GO}[\Gvf_\Gd] \|_{L^2(\GO)} \le C_2 \| \Gvf_\Gd \|_{\Hcal^*}
$$
for some positive constants $C_1$ and $C_2$, and \eqnref{reoinf2} is equivalent to
\beq\label{reoinf3}
\| \Gvf_\Gd \|_{\Hcal^*} \to \infty \quad\mbox{as } \Gd \to 0.
\eeq
It follows from \eqnref{specGvf} that
\beq\label{GvfGd}
\| \Gvf_\Gd \|_{\Hcal^*}^2 =\int_{-b}^{b} \frac{1}{|\Gl-t|^2} \,d \la \p_\nu q, \Ecal_t[\p_\nu q] \ra_{\Hcal^*}.
\eeq
As $\Gd \to 0$, $\Gl$ tends to $\Gl_0$ which is defined to be
$$
\Gl_0:= \frac{\Ge_c +1}{2(\Ge_c-1)}.
$$
The formula \eqnref{GvfGd} shows that the resonance \eqnref{reoinf3} occurs only when $\Gl_0 \in [-b, b] \subset (-1/2,1/2)$. This condition can be fulfilled only if $\Ge_c <0$. So, resonance occurs only when the dielectric constant of the inclusion is negative. A material with a negative dielectric constant is called a plasmonic material.

Let
\beq\label{muqt}
\mu_q(t):= \la \p_\nu q, \Ecal_t[\p_\nu q] \ra_{\Hcal^*}
\eeq
for ease of notation ($\mu_q$ is called the spectral measure). Since $d\mu_q$ is absolutely continuous, $\mu_q' \in L^1( -b, b )$ (Lemma \ref{specmeasure}). Observe that $\Gl= \Gl_0 + i c \Gd + O(\Gd^2)$ for some positive constant $c$. So we assume from now on that
\beq
\Gl= \Gl_0+ i\Gd
\eeq
for simplicity and $\Gl_0 \in [-b, b]$. Then \eqnref{GvfGd} can be expressed as
\beq\label{GvfGd2}
\| \Gvf_\Gd \|_{\Hcal^*}^2 =\int_{-b}^{b} \frac{\mu_q'(t)}{(\Gl_0-t)^2 + \Gd^2} \,dt.
\eeq
In what follows we investigate how fast $\| \Gvf_\Gd \|_{\Hcal^*}$ tends to $\infty$ as $\Gd \to 0$ in a quantitatively precise manner. For that, we emphasize that
$$
\frac{1}{\pi} \frac{\Gd}{(\Gl_0-t)^2 + \Gd^2}
$$
is the Poisson kernel on the half space, and it holds that
\beq\label{poisson}
\lim_{\Gd \to 0} \int_{-\infty}^\infty \frac{\Gd}{(\Gl_0-t)^2 + \Gd^2} f(t)dt = \frac{\pi}{2} [f(\Gl_0-)+f(\Gl_0+)]
\eeq
for all integrable function $f$ if the limit exists. Here $f(\Gl_0-)$ indicates the limit from the left and $f(\Gl_0+)$ from the right.

%%%%%%%%%%%%%%%%%%%%%%%%%%%%%%%%%%%%%%%%%%%%%%%%%%%%%%%%
\subsection{Generalized eigenfunctions}
%%%%%%%%%%%%%%%%%%%%%%%%%%%%%%%%%%%%%%%%%%%%%%%%%%%%%%%%

For $s \ge 0$, let $\Gd_{s}$ be the Dirac mass at $s$ and define
\beq
\Bd^1_{s} := \begin{bmatrix} \sqrt{2}\; p_1(s)^{-1/2} \Gd_{s} \\ 0 \end{bmatrix} \quad\mbox{and}\quad \Bd^2_{s} := \begin{bmatrix} 0 \\ \sqrt{2}\; p_2(s)^{-1/2}(\Gd_{s}-\Gd_0) \end{bmatrix}.
\eeq
Since $\Kbb^*$ is a multiplication operator given by \eqnref{Kbbeta}, $\Bd^1_{s}$ and $\Bd^1_{-s}$ are generalized eigenfunctions of $\Kbb^*$ corresponding to the continuous spectrum $\eta(s) \in \Gs(\Kbb^*)$, and $\Bd^2_{s}$ and $\Bd^2_{-s}$ are those corresponding to $-\eta(s)$. We emphasize that $p_j(s)^{-1/2}$ appears in the definition of $\Bd_j$ because of weights in the definition of $W^*$ and $\Gd_0$ is subtracted in the definition of $\Bd^2_s$ because $f_2(0)=0$ (if $\Bf \in W^*$ is continuous). Note that
\beq\label{eqn:1230}
\la \Bf, \Bd_s^j \ra_{W^*} =\frac{1}{\sqrt{2}} p_j(s)^{1/2} f_j(s), \quad j=1,2.
\eeq
%In particular, $\la \Bd_s^i,\Bd_t^j\ra_{W^*}=\delta_{i,j}\delta_{s,t}$ for $i,j=1,2,\ s,t\geq 0$.

Let
\beq\label{eqn:psisj}
\psi_{s}^j := U^{-1}[\Bd^j_{s}], \quad j=1,2.
\eeq
These are generalized eigenfunctions of $\Kcal_{\p\GO}^*$. We emphasize that since $\Bd_{s}^j \notin W^*$, $\psi_{s}^j \notin \Hcal^*$. In fact, as we show below, $|\psi_{s}^j (x)|$ tends to $\infty$ and $\psi_{s}^j (x)$ oscillates fast as $x$ approaches to one of the intersecting points. But, since $\Gd_s \in H^{-1}(\Rbb)$, we infer that $\psi_{s}^j \in H^{-1}(\p\GO)$. Since $\Scal_{\p\GO}$ maps $H^{-1}(\p\GO)$ into $L^2(\p\GO)$, $\la \Gvf, \psi_s^j \ra_{\Hcal^*} = - \la \Gvf, \Scal_{\p\GO} [\psi_s^j] \ra$ is well defined for $\Gvf \in L^2(\p\GO)$, and we have
\beq\label{geneigen}
\la \Gvf, \Kcal_{\p\GO}^* [\psi_s^j]  \ra_{\Hcal^*} = (-1)^{j+1} \eta(s) \la \Gvf, \psi_s^j \ra_{\Hcal^*} \quad\mbox{for all } \Gvf \in L^2(\p\GO).
\eeq

One can easily see that
$$
\Fcal^{-1}\GL^{-1}[\Bd^1_{s}](\xi) = \frac{1}{2\sqrt{\pi}} \begin{bmatrix}  p_1(s)^{-1/2} e^{is \xi} \\ -p_1(s)^{-1/2} e^{is \xi} \end{bmatrix}.
$$
So, we have from \eqnref{Cinverse}
\beq\notag
\psi_{s}^1 (x) =
\begin{cases}
\ds \frac{1}{2\sqrt{\pi}}\; p_1(s)^{-1/2} h(\Psi_1(x),\Gt_0) e^{is \Psi_1(x)}, \quad &x \in \GG_+, \\[0.2cm]
\ds -\frac{1}{2\sqrt{\pi}}\; p_1(s)^{-1/2} h(\Psi_1(x),\Gt_0) e^{is \Psi_1(x)}, \quad &x \in \GG_- .
\end{cases}
\eeq
Likewise, one can see that
\beq\notag
\psi_{s}^2 (x) = \frac{1}{2\sqrt{\pi}}\; p_2(s)^{-1/2} h(\Psi_1(x),\Gt_0) \bigr[e^{is \Psi_1(x)}-1\bigr], \quad x \in \GG.
\eeq
Since $\Psi_1(x) \to \infty$ as $x$ approaches to intersecting points of two circles, we can see from these formula that $|\psi_{s}^j (x)|$ tends to $\infty$ and $\psi_{s}^j (x)$ oscillates fast near the intersecting points.

We can also compute the single layer potentials of these function. In view of \eqnref{Sbb} and \eqnref{SbbBf} we have
$$
\Scal_{\p\GO}[\psi_{s}^j] = (\GL \Fcal \Mcal_h \Ccal)^{-1} \Sbb U [\psi_{s}^j]  = -(\GL \Fcal \Mcal_h \Ccal)^{-1} [P(s) \Bd_{s}^j],
$$
and hence
\beq\label{psione}
\Scal_{\p\GO}[\psi_{s}^1] (x) =
\begin{cases}
\ds \frac{1}{2\sqrt{\pi}}\; p_1(s)^{1/2} e^{is \Psi_1(x)}, \quad &x \in \GG_+, \\[0.2cm]
\ds -\frac{1}{2\sqrt{\pi}}\; p_1(s)^{1/2} e^{is \Psi_1(x)}, \quad &x \in \GG_-,
\end{cases}
\eeq
and
\beq\label{psitwo}
\Scal_{\p\GO}[\psi_{s}^2] (x) = \frac{1}{2\sqrt{\pi}}\; p_2(s)^{1/2} \big[e^{is \Psi_1(x)}-1\big], \quad x \in \GG.
\eeq

\begin{lemma}\label{singleform}
It holds that
\begin{align}\label{spsione}
\ds\Scal_{\p\GO}[\psi_{s}^1] (z)& = \frac{1}{2\sqrt{\pi}}\; p_1(s)^{1/2} \;\frac{\sinh s\Psi_2(z)}{\sinh s \Gt_0}\; e^{is \Psi_1(z)},\\
\label{spsitwo}
\ds\Scal_{\p\GO}[\psi_{s}^2] (z)& = \frac{1}{2\sqrt{\pi}}\; p_2(s)^{1/2} \left[ \frac{\cosh s\Psi_2(z)}{\cosh s \Gt_0} e^{is \Psi_1(z)}-1 \right],
\end{align}
for $z \in \Rbb^2 \setminus \GO$.
\end{lemma}
\pf
Since $\Psi_1(z)+i\Psi_2(z)$ is analytic, the functions on the righthand sides are harmonic in $\Rbb^2 \setminus \overline{\GO}$.
We see from \eqnref{psione} and \eqnref{psitwo} that equalities in \eqnref{spsione} and \eqnref{spsitwo} hold on $\GG=\p\GO$.
Since $\Psi_1(z)+i\Psi_2(z) \to 0$ as $|z| \to \infty$, we see that functions on the righthand side of \eqnref{spsione} and \eqnref{spsitwo}  tend to $0$ as $|z| \to \infty$. On the other hand, since $\la \psi_s^j, 1 \ra_{\Hcal^*}=0$, we infer that $\Scal_{\p\GO}[\psi_{s}^j] (z) \to 0$ as $|z| \to \infty$. So, \eqnref{spsione} and \eqnref{spsitwo} for $z \in \Rbb^2 \setminus \GO$ follow from the maximum principle. \qed

\smallskip

Let $\mu_q$ be the spectral measure defined by \eqnref{muqt}. Then we have from Lemma \ref{specmeasure} that
\begin{align*}
\mu_q'(t)=
\begin{cases}
\ds \frac{1}{|\eta'(s)|} \left[ \left| \la \p_\nu q, \psi_s^1\ra_{\Hcal^*} \right|^2  + \left| \la \p_\nu q, \psi_{-s}^1\ra_{\Hcal^*} \right|^2\right], \quad& t\in (0,b],  \\[0.3cm]
\ds\frac{1}{|\eta'(s)|}  \left[ \left| \la \p_\nu q, \psi_s^2\ra_{\Hcal^*} \right|^2  + \left| \la \p_\nu q, \psi_{-s}^2\ra_{\Hcal^*} \right|^2\right], \quad& t\in [-b, 0),
\end{cases}
\end{align*}
where $s=\eta^{-1}(t)$ if $t >0$ and $s=\eta^{-1}(-t)$ if $t <0$.
We see from Lemma \ref{singleform} that $\Scal_{\p\GO}[\psi_{-s}^j]=\overline{\Scal_{\p \GO}[\psi_s^j]}$, and hence $|\la \p_\nu q,\psi_{-s}^j\ra_{\Hcal^*}|=|\la \p_\nu q,\psi_{s}^j\ra_{\Hcal^*}|^2$ for $j=1,2$.
Thus it holds that
\beq\label{eqn:1230ddt}
\mu_q'(t) =
 \begin{cases}
 \ds\frac{2g_1(t)}{|\eta'(\eta^{-1}(t))|} ,&\quad t\in (0,b],\\[0.2cm]
\ds\frac{2g_2(-t)}{|\eta'(\eta^{-1}(-t))|} ,&\quad t\in [-b, 0).
\end{cases}
\eeq
Here we set for ease of notation
\beq
g_j(t):= \left|\la \p_\nu q,\psi_{\eta^{-1}(t)}^j\ra_{\Hcal^*}\right|^2, \quad j=1,2.
\eeq
Note that $\p_\nu q$ belongs to $H^1(\p\GO)$, is smooth except at intersecting points of two circles, real valued, and satisfies
$$
\int_{\p\GO} \p_\nu q d\Gs = \int_{\GO} \GD q(x) dx =0.
$$
So, $g_j$ is continuous on $[-b,b]$ except at $0$, and so is $\mu_q'$.

%%%%%%%%%%%%%%%%%%%%%%%%%%%%%%%%%%%%%%%%%%%%%%
\subsection{Plasmonic resonance}
%%%%%%%%%%%%%%%%%%%%%%%%%%%%%%%%%%%%%%%%%%%%%%

In this subsection we show resonance occurs at the continuous spectrum in a quantitatively precise way. Since $\mu_q'$ is continuous in $(-b, 0) \cup (0, b)$, the following proposition follows immediately from \eqnref{poisson}.

\begin{prop}\label{resoest1}
Suppose $\Gl_0 \in (-b, 0) \cup (0, b)$. Let $\Gvf_\Gd$ be the solution to {\rm\eqnref{solint}} with $\Gl=\Gl_0+i\Gd$. Then it holds that
\beq\label{reslim1}
\lim_{\Gd \to 0} \Gd \| \Gvf_\Gd \|_{\Hcal^*}^2 =
\begin{cases}
\ds \frac{2\pi g_1(\Gl_0)}{|\eta'(\eta^{-1}(\Gl_0))|}  \quad  &\mbox{if } \Gl_0 \in (0,b), \\[0.3cm]
\ds \frac{2\pi g_2(\Gl_0)}{|\eta'(\eta^{-1}(-\Gl_0))|} \quad  &\mbox{if } \Gl_0 \in (-b,0).
\end{cases}
\eeq
\end{prop}

Let us consider the case when $\Gl_0=b$. We see from \eqnref{GvfGd2}, \eqnref{poisson} and \eqnref{eqn:1230ddt} that
$$
\lim_{\Gd \to 0} \Gd \| \Gvf_\Gd \|_{\Hcal^*}^2 =
\frac{\pi}{2} \mu_q'(b-) =
\lim_{t \to b-} \frac{\pi g_1(t)}{|\eta'(\eta^{-1}(t))|} .
$$
Using \eqnref{etazero} we can see that
\beq\label{endcase2}
\lim_{t \to b-} \frac{|\eta'(\eta^{-1}(t))|}{\sqrt{b - t}} = c
\eeq
for some nonzero $c$. So we infer that
$$
\lim_{\Gd \to 0} \Gd \| \Gvf_\Gd \|_{\Hcal^*}^2 = \infty
$$
as long as $g_1(b) \neq 0$. We may refine this result in the following way.

\begin{prop}\label{resoest4}
Suppose that $\Gl_0=-b$ or $b$. Let $\Gvf_\Gd$ be the solution to {\rm\eqnref{solint}} with $\Gl=\Gl_0+i\Gd$. Then it holds that
\beq\label{reslim2}
\lim_{\Gd \to 0} \Gd^{3/2} \| \Gvf_\Gd \|_{\Hcal^*}^2 =
\begin{cases}
\ds C g_1(b)  \quad & \mbox{if }\Gl_0=b,  \\
\ds C g_2(-b) \quad & \mbox{if }\Gl_0=-b
\end{cases}
\eeq
for some nonzero constant $C$.
\end{prop}
\pf
Suppose that $\Gl_0=b$. We see from \eqnref{GvfGd2} that
\begin{align*}
\Gd^{3/2} \| \Gvf_\Gd \|_{\Hcal^*}^2 & = \int_{-b}^{b} \frac{\Gd^{3/2} \mu_q'(t)}{(b - t)^2 + \Gd^2} \, dt.
\end{align*}
Since the integral over $[-b, 0]$ tend to $0$ as $\Gd \to 0$, it suffices to consider the one over $[0, b]$. We have from \eqnref{eqn:1230ddt} that
\begin{align*}
\int_{0}^{b} \frac{\Gd^{3/2} \mu_q'(t)}{(b - t)^2 + \Gd^2} \, dt = \int_{0}^{b} \frac{\Gd^{3/2}}{(b - t)^2 + \Gd^2} \frac{2 g_1(t)}{|\eta'(\eta^{-1}(t))|} \, dt .
\end{align*}
We show that there is a constant $C \neq 0$ such that
\beq\label{Cnonzero}
\lim_{\Gd \to 0} \int_{0}^{b} \frac{\Gd^{3/2}}{(b - t)^2 + \Gd^2} \frac{\pi}{|\eta'(\eta^{-1}(t))|} \, dt = C .
\eeq
In fact, one can see from \eqnref{etazero} that
$$
\frac{1}{|\eta'(\eta^{-1}(t))|} = \frac{1}{c \sqrt{b-t}} + e(t)
$$
where $e(t)$ is a continuous function on $[0, b]$ and $c$ is the number appearing in \eqnref{endcase2}. So, we have
$$
\lim_{\Gd \to 0} \int_{0}^{b} \frac{\Gd^{3/2}}{(b - t)^2 + \Gd^2} \frac{\pi}{|\eta'(\eta^{-1}(t))|} \, dt
= \lim_{\Gd \to 0} \int_{0}^{b} \frac{\Gd^{3/2}}{(b - t)^2 + \Gd^2} \frac{\pi}{c \sqrt{b-t}} \, dt.
$$
By making change of variables, $b-t=\Gd s$, we have
$$
\lim_{\Gd \to 0} \int_{0}^{b} \frac{\Gd^{3/2}}{(b - t)^2 + \Gd^2} \frac{\pi}{c \sqrt{b-t}} \, dt
= \frac{\pi}{c} \int_{0}^{\infty} \frac{1}{(s^2 + 1) \sqrt{s}} \, ds.
$$
Using the constant $C$ in \eqnref{Cnonzero}, we write
\begin{align*}
& \int_{0}^{b} \frac{\Gd^{3/2}}{(b - t)^2 + \Gd^2} \frac{\pi g_1(t)}{|\eta'(\eta^{-1}(t))|}  \, dt - C g(b) \\
&= \int_{0}^{b} \frac{\Gd^{3/2}}{(b - t)^2 + \Gd^2} \frac{\pi [g_1(t)-g_1(b)]}{|\eta'(\eta^{-1}(t))|}  \, dt
+ g(b) \left[ \int_{0}^{b} \frac{\Gd^{3/2}}{(b - t)^2 + \Gd^2} \frac{\pi}{|\eta'(\eta^{-1}(t))|} \, dt - C \right].
\end{align*}
Using continuity of $g$ and \eqnref{Cnonzero}, it can be proved by standard estimates that terms on the righthand side above
converge to $0$ as $\Gd \to 0$. So \eqnref{reslim2} for $\Gl_0=b$ is proved. The case when $\Gl_0=-b$ can be handled similarly. \qed

Observe that the case when $\Gl_0=0$ is missing in Proposition's of this subsection. We deal with this case when the source function is a polarized dipole in the next section.
Results of this subsection show that if
\beq\label{nonzerocon}
g_j(\Gl_0) \neq 0,
\eeq
then resonance occurs at the order of $\Gd^{-1/2}$ or $\Gd^{-3/4}$, namely $\| \Gvf_\Gd \|_{\Hcal^*} \to \infty$ at the order of $\Gd^{-1/2}$ or $\Gd^{-3/4}$ as $\Gd \to 0$. We will show in the next section that resonance occurs at $\Gl_0=0$ at the rate of $\Gd^{-\Ga}$ for any $\Ga<1$ where $\Ga$ varies depending on the location of the dipole source. We will also show that \eqnref{nonzerocon} is fulfilled with polarized dipole sources.

It is worthwhile to mention that this order of resonance (at the continuous spectrum) is in contrast with that at the discrete spectrum where the order is $\Gd^{-1}$ (see \cite{Ando}). The following lemma shows that at the absolutely continuous spectrum (on any domain with Lipschitz boundary) resonance never occur at the order of $\Gd^{-1}$.

\begin{prop}\label{never}
Let $\GO$ be a bounded domain with the Lipschitz boundary. Suppose that the NP operator on $\GO$ has only absolutely continuous spectrum. Let $\Gvf_\Gd$ be the solution to {\rm\eqnref{solint}} with $\Gl=\Gl_0+i\Gd$. Then, it holds that
\beq
\lim_{\Gd \to 0} \Gd \| \Gvf_\Gd \|_{\Hcal^*} = 0.
\eeq
\end{prop}
\pf
Let $\mu$ be the spectral measure (as defined in \eqnref{muqt}) for the NP operator on $\GO$. Since only absolutely continuous spectrum is nonempty, $\mu' \in L^1(\Rbb^1)$, and
it holds that
\beq
\| \Gvf_\Gd \|_{\Hcal^*}^2 = \int_{-\infty}^{\infty} \frac{\mu'(t)}{(\Gl_0 - t)^2 + \Gd^2} \, dt .
\eeq
For simplicity, we assume $\Gl_0=0$. We decompose $\Gd^2 \| \Gvf_\Gd \|_{\Hcal^*}^2$ as
$$
\Gd^2 \| \Gvf_\Gd \|_{\Hcal^*}^2 = \int_{-\infty}^{\infty} \frac{\Gd^{2}}{t^2 + \Gd^2} |\mu'(t)| \, dt  = \int_{|t| \le \Gd} + \sum_{k=1}^\infty \int_{2^{k-1} \Gd < |t| \le 2^{k}\Gd}= :I_0 + \sum_{k=1}^\infty  I_k.
$$
Since $\mu'$ is integrable, we have $I_0 \to 0$ as $\Gd \to 0$. For $k > 0$, we have
$$
I_k = \int_{2^{k-1} \Gd < |t| \le 2^{k}\Gd} \frac{\Gd^{2}}{t^2 + \Gd^2} \mu'(t) \, dt \le \frac{1}{1+ 2^{2(k-1)}} \int_{2^{k-1} \Gd < |t| \le 2^{k}\Gd}  \mu'(t) \, dt.
$$
Observe that
$$
\int_{2^{k-1} \Gd < |t| \le 2^{k}\Gd}  \mu'(t) \, dt \to 0
$$
as $\Gd \to 0$. So by Lebesgue dominated convergence theorem, we have
$$
\sum_{k=1}^\infty \frac{1}{1+ 2^{2(k-1)}} \int_{2^{k-1} \Gd < |t| \le 2^{k}\Gd}  \mu'(t) \, dt \to 0
$$
as $\Gd \to 0$. This completes the proof. \qed

%%%%%%%%%%%%%%%%%%%%%%%%%%%%%%%%%%%%%%%%%%%%%%%%%%%%%%%%%%%%
\section{Resonance by dipole sources}
%%%%%%%%%%%%%%%%%%%%%%%%%%%%%%%%%%%%%%%%%%%%%%%%%%%%%%%%%%%%

Suppose that the source $f$ is given by
\beq\label{dipole1}
f(x) =a \cdot \nabla \Gd_z(x)
\eeq
where $a$ is a constant vector and $z$ is a point outside $\GO$. In this case, we have
\beq
q(x) = a \cdot \nabla_x \GG(x-z)= -a \cdot \nabla_z \GG(x-z).
\eeq
So we have
\beq\label{lapnu}
\la \p_\nu q, \psi_{s}^j \ra_{\Hcal^*} = a \cdot \nabla_z \int_{\p\GO} \pd{}{\nu_y} \GG(z-y) \Scal_{\p\GO}[\psi_{s}^j](y) \, d\Gs(y).
\eeq
Note that
\begin{align*}
\int_{\p\GO} \pd{}{\nu_y} \GG(z-y) \Scal_{\p\GO}[\psi_{s}^j](y) \, d\Gs(y) = \Dcal_{\p\GO} \Scal_{\p\GO}[\psi_{s}^j](z)
\end{align*}
Recall that the double layer potential $\Dcal_{\p\GO}$ satisfies
$$
\Dcal_{\p\GO}[\Gvf] |_{+} = (-\frac{1}{2}I + \Kcal_{\p\GO})[\Gvf]
$$
for all $\Gvf \in L^2(\p\GO)$. So, we have from \eqnref{geneigen}
\begin{align*}
\la \Gvf, \Dcal_{\p\GO}\Scal_{\p\GO}[\psi_{s}^j] |_{+} \ra &= \big\la \Gvf, \big(-\frac{1}{2}I + \Kcal_{\p\GO}\big)\Scal_{\p\GO}[\psi_{s}^j] \big\ra \\
&= \big\la \big(-\frac{1}{2}I + \Kcal_{\p\GO}^*\big)[\Gvf], \Scal_{\p\GO}[\psi_{s}^j] \big\ra \\
&= \big\la \big(\frac{1}{2}I - \Kcal_{\p\GO}^*\big)[\Gvf], \psi_{s}^j \big\ra_{\Hcal^*} \\
&= \big( \frac{1}{2} - (-1)^{j+1} \eta(s) \big) \la \Gvf, \psi_{s}^j \ra_{\Hcal^*}  \\
&= \big( -\frac{1}{2} + (-1)^{j+1} \eta(s) \big) \la \Gvf, \Scal_{\p\GO}[\psi_{s}^j] \ra.
\end{align*}
So we have
$$
\Dcal_{\p\GO}\Scal_{\p\GO}[\psi_{s}^j] = \big( -\frac{1}{2} + (-1)^{j+1} \eta(s) \big) \Scal_{\p\GO}[\psi_{s}^j] \quad\mbox{on } \p\GO.
$$
Since $\Dcal_{\p\GO}\Scal_{\p\GO}[\psi_{s}^j](z)$ and $\Scal_{\p\GO}[\psi_{s}^j](z)$ are harmonic in $\Rbb^2 \setminus \overline{\GO}$ and decay to $0$ as $|z| \to \infty$, we infer that
\begin{align*}
\int_{\p\GO} \pd{}{\nu_y} \GG(z-y) \Scal_{\p\GO}[\psi_{s}^j](y) \, d\Gs(y) = \big( -\frac{1}{2} + (-1)^{j+1} \eta(s) \big) \Scal_{\p\GO}[\psi_{s}^j](z).
\end{align*}
We then see from \eqnref{lapnu} that
\beq\label{gjeta}
g_j(\eta(s)) = \left| \la \p_\nu q, \psi_{s}^j \ra_{\Hcal^*} \right|^2 = \left| -\frac{1}{2} + (-1)^{j+1} \eta(s) \right|^2 \left|
a \cdot \nabla \Scal_{\p\GO}[\psi_{s}^j](z) \right|^2.
\eeq

Choose $z \in \Rbb^2 \setminus \overline{\GO}$ so that
\beq\label{dipole3}
\Psi_2(z) \neq 0,
\eeq
which is equivalent to that $z$ is not on the real axis.
We see from Lemma \ref{singleform} that
\begin{align}
& \left| a \cdot \nabla \Scal_{\p\GO}[\psi_{s}^1] (z) \right|^2 \nonumber \\
& = \frac{1}{4\pi} \frac{s^2 p_1(s)}{\sinh^2 s \Gt_0} \Big[ (a \cdot \nabla \Psi_2(z))^2 \cosh^2 s \Psi_2(z) +  (a \cdot \nabla \Psi_1(z))^2 \sinh^2 s \Psi_2(z) \Big]. \label{dipole4}
\end{align}
Since $\Psi_1(z) + i \Psi_2(z)$ is conformal, we have
$$
(a \cdot \nabla \Psi_1(z))^2 + (a \cdot \nabla \Psi_2(z))^2 \neq 0
$$
for any $a$. We choose $a$ so that
\beq\label{dipole5}
a \cdot \nabla \Psi_1(z) \neq 0 \quad\mbox{and}\quad a \cdot \nabla \Psi_2(z) \neq 0.
\eeq
Then, $| a \cdot \nabla \Scal_{\p\GO}[\psi_{s}^1] (z)| \neq 0$ for all $s$, and hence $g_1(\eta(s)) \neq 0$, or \eqnref{nonzerocon} is satisfied.
Similarly, we obtain
\begin{align*}
& \left| a \cdot \nabla \Scal_{\p\GO}[\psi_{s}^2] (z) \right|^2 \\
& = \frac{1}{4\pi} \frac{s^2 p_2(s)}{\cosh^2 s \Gt_0} \Big[ (a \cdot \nabla \Psi_2(z))^2 \sinh^2 s \Psi_2(z) +  (a \cdot \nabla \Psi_1(z))^2 \cosh^2 s \Psi_2(z) \Big],
\end{align*}
and hence $g_2(\eta(s)) \neq 0$.

Let us now consider the resonance at $\Gl_0 =0$. To do that we need to look into the behavior of
$g_j(t)/|\eta'(\eta^{-1}(t))|$ as $t \to 0+$. We see from \eqnref{gjeta} and \eqnref{dipole4} that
$$
g_1(\eta(s)) = C_1 s e^{2s (-\Gt_0 + |\Psi_2(z)|)} \big[1+ O(e^{-2 |\Psi_2(z)| s})\big] \quad\mbox{as } s \to \infty
$$
for some constant $C_1 \neq 0$.
We also see from the definition \eqnref{etadef} of $\eta$ that
$$
\eta(s) =  C_2 e^{-2\Gt_0 s} \big[1+ O(e^{-2(\pi-2\Gt_0)s})\big]
$$
and
$$
\eta'(s) =  -2\Gt_0 C_2 e^{-2\Gt_0 s} \bigr[1+ O(e^{-2(\pi-2\Gt_0)s})\bigr]
$$
as $s \to \infty$ for some constant $C_2 \neq 0$.
Therefore we have
$$
\frac{\pi g_1(\eta(s))}{2|\eta'(s)|} = C_3 s e^{2s |\Psi_2(z)|} \big[1+ O(e^{-2 \Ge s})\big]
$$
for some constant $C_3 \neq 0$ and $\Ge>0$, which implies that
\beq\label{dipole7}
\frac{\pi g_1(t)}{2|\eta'(\eta^{-1}(t))|}  = C_3  |\log t| t^{-|\Psi_2(z)|/\Gt_0} + e_1(t)
\eeq
as $t \to 0+$, where the error term $e_1(t)$ satisfies
\beq\label{bone}
|e_1(t)| \le C |\log t| t^{-|\Psi_2(z)|/\Gt_0 + \Ge}.
\eeq
for some constant $C$. Likewise we obtain
\beq\notag
\frac{\pi g_2(-t)}{2|\eta'(\eta^{-1}(-t))|} = C_4  |\log t| t^{-|\Psi_2(z)|/\Gt_0} + e_2(t)
\eeq
as $t \to 0+$ for some $C_4 \neq 0$, where
$e_2(t)$ satisfies
\beq\label{btwo}
|e_2(t)| \le C |\log t| t^{-|\Psi_2(z)|/\Gt_0 + \Ge}.
\eeq

\begin{theorem}\label{thm:dipole2}
Suppose that the source $f$ is given by {\rm\eqnref{dipole1}} with $a$ and $z$ chosen to satisfy {\rm\eqnref{dipole3}} and {\rm\eqnref{dipole5}}.
Let $\Gvf_\Gd$ be the solution to {\rm\eqnref{solint}} with $\Gl=i\Gd$.  Then there is a constant $C>0$ such that
\beq\label{resoest10}
|\log\Gd|^{-1} \Gd^{1+|\Psi_2(z)|/\Gt_0} \| \Gvf_\Gd \|_{\Hcal^*}^2 \to C
\eeq
as $\Gd \to 0$
\end{theorem}
\pf
We have from \eqnref{GvfGd}
\begin{align*}
\| \Gvf_\Gd \|_{\Hcal^*}^2 & = \int_{-b}^{b} \frac{\mu_q'(t)}{t^2 + \Gd^2} \, dt \\
&= \int_{0}^{b} \frac{1}{t^2 + \Gd^2} \frac{\pi g_1(t)}{2|\eta'(\eta^{-1}(t))|} \, dt + \int_{0}^{b} \frac{1}{t^2 + \Gd^2} \frac{\pi g_2(t)}{2|\eta'(\eta^{-1}(t))|}  \, dt \\
&=: I_1(\Gd) + I_2(\Gd).
\end{align*}
We show below that
\beq\label{Ij}
|\log\Gd|^{-1} \Gd^{1+|\Psi_2(z)|/\Gt_0} I_j(\Gd) \to C_j \neq 0
\eeq
as $\Gd \to 0$ for $j=1,2$. Then \eqnref{resoest10} follows immediately.

Let us prove \eqnref{Ij} for $j=1$. The case for $j=2$ can be handled with in the same way.
Let $e_1$ be the function appeared in \eqnref{dipole7}. Then one can see from \eqnref{bone}, by the scaling $t =\Gd s$, that
\beq\label{dipolezero1}
\lim_{\Gd \to 0}|\log\Gd|^{-1} \Gd^{1+|\Psi_2(z)|/\Gt_0} \int_{0}^{b} \frac{e_1(t)}{t^2 + \Gd^2}  \, dt = 0.
\eeq
On the other hand we see by the same scaling that
\begin{align*}
|\log\Gd|^{-1} \Gd^{1+|\Psi_2(z)|/\Gt_0}\int_{0}^{b} \frac{1}{t^2 + \Gd^2} |\log t| t^{-|\Psi_2(z)|/\Gt_0} \, dt &= \int_{0}^{b/\Gd} \frac{1}{s^2 + 1} \left| \frac{\log s}{\log\Gd} +1 \right| s^{-|\Psi_2(z)|/\Gt_0} \, ds \\
&= \int_{0}^{\Gd} + \int_{\Gd}^{b/\Gd} =: II_1(\Gd)+ II_2(\Gd).
\end{align*}
If $s >\Gd$, then $\log s / \log\Gd \to 0$ as $\Gd \to 0$ for all $s$. So by Lebesgue dominated convergence theorem we have
\beq\label{dipolezero2}
II_2(\Gd) \to \int_{0}^{\infty} \frac{1}{s^2 + 1} s^{-|\Psi_2(z)|/\Gt_0} \, ds \neq 0
\eeq
as $\Gd \to 0$. If $0< s \le \Gd$, then $|\log s / \log\Gd +1| \le 2|\log s|$ if $\Gd$ is sufficiently small, and hence
$$
II_1(\Gd) \le 2 \int_{0}^{\Gd} \frac{1}{s^2 + 1} |\log s| s^{-|\Psi_2(z)|/\Gt_0} \, ds.
$$
Since $|\Psi_2(z)|<\Gt_0$, the integrand on the righthand side is integrable, and it follows that
\beq\label{dipolezero3}
II_1(\Gd) \to 0
\eeq
as $\Gd \to 0$. By combining \eqnref{dipolezero2} and \eqnref{dipolezero3} we infer that
$$
|\log\Gd|^{-1} \Gd^{1+|\Psi_2(z)|/\Gt_0} \int_{0}^{b} \frac{1}{t^2 + \Gd^2} |\log t| t^{-|\Psi_2(z)|/\Gt_0} \, dt \to \int_{0}^{\infty} \frac{1}{t^2 + 1} t^{-|\Psi_2(z)|/\Gt_0} \, dt
$$
as $\Gd \to 0$, which together with \eqnref{dipolezero1} implies \eqnref{Ij} for $j=1$. This completes the proof. \qed

Estimate \eqnref{resoest10} shows that resonance occurs at $\Gl_0=0$ at the rate of $|\log\Gd|^{1/2} \Gd^{-1/2-|\Psi_2(z)|/2\Gt_0}$. So, the stronger resonance occurs if $|\Psi_2(z)|/\Gt_0$ is closer to 1, in other words, if $z$ is closer to $\p\GO$. But since $|\Psi_2(z)| <\Gt_0$, the resonance at the rate of $\Gd^{-1}$ is never achieved as we showed in Proposition \ref{never}.

%%%%%%%%%%%%%%%%%%%%%%%%%%%%%%%%%%


\begin{thebibliography}{09}

\bibitem{ACKLM} H. Ammari, G. Ciraolo, H. Kang, H. Lee and G.W. Milton, Spectral theory of a Neumann-Poincar\'e-type operator and analysis of cloaking due to anomalous localized resonance, Arch. Ration. Mech. An. 208 (2013), 667--692.

\bibitem{ACKLY} H. Ammari, G. Ciraolo, H. Kang, H. Lee and K. Yun, Spectral analysis of the Neumann-Poincar\'e operator and characterization of the stress concentration in anti-plane elasticity, Arch. Rati. Mech. Anal. 208 (2013), 275--304.

\bibitem{AmKa07Book2} {H. Ammari and H. Kang}, {\sl Polarization and moment tensors},
{Applied Mathematical Sciences}, 162, Springer, New York. 2007.

\bibitem{Ando} K. Ando and H. Kang, Analysis of plasmonic resonance on smooth domains using spectral properties of the Neumann-Poincar\'e operator, submitted, arXiv 1412.6250.

\bibitem{BT} E. Bonnetier and F. Triki, Pointwise bounds on the gradient and the spectrum of the
Neumann-Poincar\'e operator: The case of 2 discs, Contemp. Math. 577 (2012), 79--90.

\bibitem{BT2} E. Bonnetier and F. Triki, On the spectrum of Poincar\'e variational problem for two close-to-touching inclusions in 2D, Arch. Ration. Mech. An. 209 (2013), 541--567.

\bibitem {Fo95} {G. B. Folland},
{\sl Introduction to partial differential equations}, 2nd Ed., Princeton Univ. Press, Princeton, 1995.

\bibitem{Grieser} D. Grieser, Plasmonic eigenvalue problem, Rev. Math. Phys. 26 (2014), 1450005.

\bibitem{Kang}
   {H. Kang},
    {Layer potential approaches to interface problems},
     {\sl Inverse problems and imaging, Panoramas et Syntheses},
 {Societe Mathematique de France},
 {to appear}.

\bibitem{KKLS} H. Kang, K. Kim, H. Lee and J. Shin, Spectral properties of the Neumann-Poincar\'e operator and uniformity of estimates for the conductivity equation with complex coefficients, submitted, arXiv 1406.3873.

\bibitem{KS2000}
H. Kang and J-K. Seo, {Recent progress in the inverse conductivity problem with
single measurement}, in Inverse Problems and Related Fields, CRC
Press, Boca Raton, FL, 2000, 69--80.

\bibitem{Ke29} {O. D. Kellogg},
    {\sl Foundations of potential theory}, Dover, New York, 1953 (Reprint from the first edition of {\sl Die Grundlehren der
              Mathematischen Wissenschaften}, Band 31,
 {Springer-Verlag, Berlin-New York}, 1929).

\bibitem{kenig} C. Kenig, Harmonic analysis techniques for second order elliptic boundary value problems, CBMS series 83, Amer. Math. Soc. Providence, 1991.

\bibitem{KhPuSh07} {D. Khavinson, M. Putinar,  and H. S. Shapiro},
    {Poincar\'e's variational problem in potential theory},
  {Arch. Ration. Mech. An.} 185 (2007), 143--184.

\bibitem{pendry} D. Y. Lei, A. Aubry, Y. Luo, S. A. Maier, and J. B. Pendry, Plasmonic interaction between overlapping nanowires, ACS Nano 5(1) (2011), 597-–607.

\bibitem {Li01} {M. Lim},
     {Symmetry of a boundary integral operator and a
              characterization of a ball},
  {Illinois J. Math.} 45 (2001), 537--543.

\bibitem{LY1} M. Lim and S. Yu, Asymptotics of the solution to the conductivity equation in
the presence of adjacent circular inclusions with finite
conductivities, J. Math. Anal. Appl. 421 (2015), 131-156.

\bibitem{MM} R. C. McPhedran and G. W. Milton, Transport Properties of Touching Cylinder Pairs and of the Square Array of Touching Cylinders, Proc. R. Soc. Lond. A 411 (1987), 313-326.

\bibitem{PP} K.-M. Perfekt and M. Putinar, Spectral bounds for the Neumann-Poincar\'e operator on planar domains with corners,  J. Anal. Math. 124 (2014), 39-57.

\bibitem{RS} M. Reed and B. Simon, {\sl Methods of modern mathematical physics. I. Functional analysis}, Revised and enlarged edition, Academic Press, New York, 1980.

\bibitem{Yo} K. Yosida, {\sl Functional Analysis}, 4th Ed., Springer, Berlin, 1974.


\end{thebibliography}
\end{document}